\documentclass[11pt]{article}
\usepackage{graphicx}
\usepackage{color}
\usepackage{amssymb}
\usepackage{amsmath}
\usepackage{epic}
\usepackage{amscd}
\usepackage{amsthm}
\usepackage{graphicx}
\usepackage{multirow}
\usepackage{float}

\newtheorem{alg}{\bf Algorithm}
\newtheorem{theorem}{\bf Theorem}

\newtheorem{defi}{\bf Definition}

\title{Implicitly Restarted Generalized Second-order Arnoldi Type Algorithms
for the Quadratic Eigenvalue Problem\footnote{Supported
by National Basic Research Program of China 2011CB302400 and the
National Science Foundation of China (Nos. 11071140, 11201020 and 11371219).}}
\author{Zhongxiao Jia\thanks{Department of Mathematical Sciences,
Tsinghua University, Beijing 100084, People's Republic of China,
{\sf jiazx@tsinghua.edu.cn}} \and Yuquan Sun\thanks{Corresponding author.
LMIB \& School of Mathematics and Systems Science, BeiHang University, Beijing
100191, People's Republic of China, {\sf sunyq@buaa.edu.cn}}}
\date{}

\begin{document}
\maketitle
\setcounter{page}{1}
\setcounter{equation}{0}
\parbox{12cm}

\begin{abstract}
We investigate the generalized second-order Arnoldi (GSOAR) method,
a generalization of the SOAR method proposed by Bai and Su
[{\em SIAM J. Matrix Anal. Appl.}, 26 (2005): 640--659.],
and the Refined GSOAR (RGSOAR) method for the quadratic
eigenvalue problem (QEP). The two methods use the GSOAR procedure to
generate an orthonormal basis of a given generalized second-order
Krylov subspace, and with such basis they project the QEP onto the
subspace and compute the Ritz pairs and the refined Ritz pairs, respectively.
We develop implicitly restarted GSOAR and RGSOAR algorithms, in which
we propose certain exact and refined shifts for respective use within
the two algorithms. Numerical experiments on real-world problems illustrate
the efficiency of the restarted
algorithms and the superiority of the restarted RGSOAR to
the restarted GSOAR. The experiments also demonstrate
that both IGSOAR and IRGSOAR generally perform much better than
the implicitly restarted Arnoldi method applied to the corresponding
linearization problems, in terms of the accuracy and the computational
efficiency.


\end{abstract}

\section{Introduction} {\normalsize\ \ \ \ \

Consider the large QEP
\begin{equation}
Q(\lambda)x=(\lambda^2 M + \lambda C +K)x=0 \label{05-qep}
\end{equation}
with $\|x\|=1$,
where $M,\ C,\ K$ are $n\times n$ matrices with $M$ nonsingular and
$\|\cdot\|$ is the 2-norm of a vector or matrix. Such QEP arises in a
wide variety of scientific and engineering applications
\cite{Betcke,tm}. One is often interested in a few largest eigenvalues
in magnitude or a few eigenvalues nearest to a target $\sigma$ in
the complex plane. One of the commonly used approaches is to linearize the
QEP and then solve the linearized problem.
There are a number of linearizations available \cite{tm}, of which
a commonly used one is to transform (\ref{05-qep}) to the
generalized eigenvalue problem
\begin{equation}
\left [ \begin{array}{cc} -C & -K \\I & 0\end{array} \right ] \left
[\begin{array}{c} \lambda x \\x \end{array} \right ]= \lambda  \left
[ \begin{array}{cc} M & 0 \\0 & I \end{array} \right ] \left [
\begin{array}{c} \lambda x \\x \end{array} \right ], \label{5xian}
\end{equation}
which is equivalent to the standard linear eigenvalue problem
\begin{equation}
\left [ \begin{array}{cc} A &B \\I&0\end{array} \right ] \left [
\begin{array}{c} \lambda x \\x \end{array} \right ]=
\lambda \left [
\begin{array}{c} \lambda x \\x \end{array} \right ], \label{5biao}
\end{equation}
where $A=-M^{-1}C, B=-M^{-1}K$. Clearly, (\ref{5biao}) corresponds to
the monic QEP
\begin{equation}
 (\lambda^2 I- \lambda A -B)x=0. \label{monic}
\end{equation}

The mathematical theory on (\ref{5xian}) and (\ref{5biao}) has been
well established and a number of numerical methods have been
available for solving them \cite{bai,golub,saad1,stewart,vorst}.
One of the drawbacks via linearizations is that general numerical methods
do not take the structures of (\ref{5xian}) and (\ref{5biao}) into account, making
computations expensive and the approximate eigenpairs possibly lose their
physical structures.

To improve the computational efficiency of the Arnoldi method that is directly
applied to
some linearization problem of QEP (\ref{05-qep}), Meerbergen \cite{meerbergen}
proposes a quadratic Arnoldi (Q-Arnoldi) method, which exploits the structure
of the linearization problem to reduce the memory requirements by about a half and
can compute a partial Schur form of the linearization problem
with respect to the structure of the Schur vectors. He shows that the Q-Arnoldi
method can be implicitly restarted. Some similar methods have proposed very recently
in \cite{Kressner, zhangs}. All of them have
are special Arnoldi methods applied to
certain linearization problems of QEP (\ref{05-qep}), that is, each of them
projects the corresponding linearization problem rather than QEP (\ref{05-qep}) or
(\ref{monic}) onto some Krylov subspace whose
orthonormal basis is generated efficiently by a special Arnoldi process.
For these methods, the implicit restarting technique \cite{so}
is easily applied.

In this paper, we are interested in projection methods that work on QEP (\ref{05-qep})
directly other than its linearizations, and such methods preserve some important
structures of it.
The second-order Arnoldi (SOAR) method proposed by Bai and Su \cite{bs} falls into
this category and is a Rayleigh--Ritz method. They propose a SOAR procedure
that computes an orthonormal basis of a second-order Krylov
subspace generated by the matrices $A$ and $B$ simultaneously.
The SOAR method then projects (\ref{05-qep}) onto this subspace
and computes the Ritz pairs to approximate the desired
eigenpairs of (\ref{05-qep}). A unified and general convergence
theory has recently been established in \cite{huangjialin} for the
Rayleigh--Ritz method and the refined Rayleigh--Ritz method for the QEP,
generalizing some of the known results on the
Rayleigh--Ritz method for the linear eigenvalue problem \cite{stewart}.
It is proved in \cite{huangjialin} that for a sequence of projection subspaces
containing increasingly accurate approximations to a desired eigenvector
there is a Ritz value that converges to the desired eigenvalue unconditionally
while the corresponding Ritz vector converges conditionally and may fail to
converge. Alternatively, we can compute a refined Ritz vector whose
unconditional convergence is guaranteed.

In the spirit of the Hessenberg-triangular
decomposition of a matrix pencil, which reduces to the Hessenberg decomposition of
a single matrix, Huang {\em et al}. \cite{wqhuang} propose a semiorthogonal
generalized Arnoldi (SGA) procedure for the matrix pencil resulting
from some linearization of the QEP.
The SGA method first generates an SGA decomposition and then computes
the Rayleigh--Ritz approximations
of QEP (\ref{05-qep}) with respect to the subspace
defined by an orthonormal basis generated by the SGA decomposition.
To overcome the possible non-convergence of Ritz vectors obtained by the
SGA method, they apply the refined projection principle \cite{J44} to
propose a refined SGA (RSGA) method that computes better refined Ritz vectors.
On the basis of implicitly shifted QZ iterations, they have
developed the implicitly restarted SGA and RSGA algorithms, abbreviated
as IRSGA and IRRSGA, with certain exact shifts and refined shifts
suggested, respectively.

One disadvantage of SOAR is that the implicit restarting technique is not
directly applicable. In order to make implicit restarting applicable,
Otto \cite{oc} proposes a
modified SOAR procedure that replaces the original special starting vector
by a general one. Under the assumption that there is no deflation in the
modified SOAR procedure, implicit restarting is directly
adapted to this procedure. However, it is hard to
interpret and understand the modified SOAR method. This is
unlike the SOAR method, whose convergence is related to
the Arnoldi method for the linear eigenvalue problem.
Wei {\em et al}. \cite{bao,zhou} make a similar modification
and propose a generalized second-order Arnoldi (GSOAR) method and their refined
variants, for solving the QEP and higher degree polynomial eigenvalue problems.
Based on the explicit restarting scheme of the
Arnoldi algorithm for the linear eigenvalue problem, Wei {\em et al}. \cite{bao,zhou}
have developed explicitly restarted generalized Krylov subspace algorithms.
Deflation and breakdown may take place in the SOAR and modified SOAR procedures,
but they have completely different consequences \cite{bs,oc}, where it is proved
that the SOAR method will find some exact eigenpairs of QEP (\ref{05-qep})
if breakdown occurs but no eigenpair is found generally when deflation takes place.
A remedy strategy is given in \cite{bs} to treat the deflation so as to continue
the SOAR procedure. Similarly, deflation may occur in the GSOAR procedure, but it
is not mentioned in \cite{bao,zhou}.

Similar to the modified SOAR procedure, implicit restarting is directly
adapted to the GSOAR procedure, but it is useable only conditionally and
requires that no deflation occur in implicit restarts. Once deflation
takes place, implicit restarting fails completely. Therefore, one must
cure deflations, so that implicit restarting can be applied unconditionally.
For the success and overall performance of implicitly restarted GSOAR type
algorithms, just as the mechanism for those implicitly restarted Krylov
subspace algorithms for the linear eigenvalue problem and SVD problems
\cite{hochsten02,J50,jia02,jianiu03,jianiu10,kokio04},
it turns out that a proper selection of the shifts is crucial.
Otto \cite{oc} has proposed certain exact shifts for his implicitly
restarted modified SOAR algorithm, but they could cause convergence problems since
some important aspects on the QEP are ignored when determining the shifts.

In this paper, we are concerned with the GSOAR and RGSOAR methods and their implicit
restarting. We will explore more properties and features of them, and consider
the efficient and reliable computation of refined Ritz vectors.
Particularly, we show that there is a close relationship between the subspace
generated by the GSOAR procedure and a standard Krylov subspace.
With help of this result, we can interpret the convergence of the GSOAR type
methods. Our main concern is a reasonable selection of the shifts when
implicitly restarting the GSOAR and RGSOAR algorithms. We advance
certain exact shifts, different from those in \cite{oc}, and refined
shifts for respective use within the implicitly restarted GSOAR
and RGSOAR algorithms. The refined shifts are based on the refined Ritz vectors
and theoretically better than the exact shifts. Unlike the implicitly restarted
algorithms for the linear eigenvalue problem, both exact and refined shift
candidates are now more than the shifts allowed.
We show how to reasonably select the desired shifts
among them. We present an efficient algorithm to compute
the exact and refined shift candidates reliably. In addition,
we propose an effective approach to cure deflations
in implicit restarts, so that implicit restarting is useable unconditionally.

The rest of this paper is organized as follows. In Section~\ref{s2} we
review the SOAR and GSOAR procedures, present some properties of them,
and describe the SOAR and GSOAR methods. In Section~\ref{s3}, we describe the RGSOAR
method and discuss some practical issues of it. In Section~\ref{s4},
we develop implicitly restarted GSOAR and RGSOAR algorithms with
the exact and refined shifts suggested. We present
an effective approach to treat deflations in implicit restarts.
In Section~\ref{s5}, we report numerical experiments
to illustrate the efficiency of the restarted algorithms and the
superiority of the refined algorithm. We also compare our algorithms with
IRSGA and IRRSGA \cite{wqhuang}, demonstrating that ours perform better.
More importantly, we compare our algorithms with the Matlab function
{\sf eigs}, the implicitly restarted Arnoldi method applied to a commonly
used linearization problem, showing that ours generally
have sharp superiority to {\sf eigs}
in terms of the accuracy and the computational efficiency.
Finally, we conclude the paper in Section 6.

Throughout the paper, we denote by $\|\cdot\|$ the spectral norm of
a matrix and the 2-norm of a vector, by $I$
the identity matrix with the order clear from the context, by the
superscripts $T$ and $*$ the transpose and conjugate transpose of a vector or
matrix, by ${\cal C}^k$ the complex vector space of dimension
$k$ and by ${\cal C}^{(k+1)\times k}$ the set of $(k+1)\times k$
matrices. We denote by $\sigma_{\min}(F)$ the smallest singular value
of a matrix $F$ and by the Matlab notation $A(i:j,k:l)$
the submatrix consisting of rows $i$ to $j$ and
columns $k$ to $l$ of $A$.

\section{The SOAR and GSOAR methods}\label{s2}

Bai and Su \cite{bs} introduce the following concepts.

\begin{defi} \label{def6.1}
Let $A,\ B$ be matrices of order $n$ and for the vector $u\neq 0$, and define
$$
\begin{array}{ccl}
r_0&=&u,\\
r_1&=&Ar_0,\\
r_j&=&Ar_{j-1}+Br_{j-2} \hspace{0.5cm} \mbox{for } j\geq 2.
\end{array}
$$
Then
$r_0,r_1,r_2,\ldots, r_{k-1}$
is called a second-order Krylov sequence based on $A, B$ and $u$, and
$
\mathcal{G}_k(A,B;u)=span\{r_0,r_1,r_2,\ldots, r_{k-1}\}
$
a $k$-th second-order Krylov subspace.
\end{defi}

Note that (\ref{5biao}) is a linearization of (\ref{monic}). Define
the matrix
\begin{equation}
H=\left[\begin{array}{cc}A&B\\I&0\end{array} \right] \label{line-1}
\end{equation}
of order $2n$. For a $2n$-dimensional starting vector $v$, we can
generate a Krylov subspace
$\mathcal{K}_k(H, v)=span\{v,  Hv,  H^2v,  \ldots, H^{k-1}v\}$.
Particularly, if we choose $v=[u^T,\ 0]^T$, we have
\begin{equation}
\left[\begin{array}{c}r_j\\r_{j-1} \end{array} \right]=H^jv, \ j\geq
 0 \mbox{ with $r_{-1}=0$}. \label{linsub}
\end{equation}
We observe the fundamental relation
\begin{equation}
{\cal K}_k(H,v)\subseteq {\cal G}_k^2(A,B;u), \label{gk1}
\end{equation}
where ${\cal G}_k^2(A,B;u)$ is the subspace generated by the vector set
$$
\left\{\left[\begin{array}{c}
r_0\\
0\\
\end{array}
\right],
\left[\begin{array}{c}
r_1\\
0\\
\end{array}
\right],
\ldots,
\left[\begin{array}{c}
r_{k-1}\\
0\\
\end{array}
\right],
\left[\begin{array}{c}
0\\
r_0\\
\end{array}
\right],
\left[\begin{array}{c}
0\\
r_1\\
\end{array}
\right],
\ldots,
\left[\begin{array}{c}
0\\
r_{k-1}\\
\end{array}
\right]
\right\}.
$$
Due to the equivalence of QEP~(\ref{monic}) and the eigenproblem of $H$, relation
(\ref{gk1}) shows that if the eigenvector $[\lambda x^T,x^T]^T$ is
contained in ${\cal K}_k(H,[u^T,0]^T)$, then the eigenvector $x$ of
QEP~(\ref{monic}) is contained in ${\cal G}_k(A,B;u)$. By continuity,
if there is a good
approximation to $[\lambda x^T,x^T]^T$ in ${\cal K}_k(H,[u^T,0]^T)$, then there
must be a good approximation to $x$ contained in ${\cal G}_k(A,B;u)$.

Bai and Su \cite{bs} propose the following procedure for computing an orthonormal
basis $\{q_j\}_{j=1}^k$ of $\mathcal{G}_k(A,B;u)$ and an auxiliary vector
sequence $\{p_j\}$ generating $\mathcal{G}_{k-1}(A,B;u)$.

\begin{alg} {\bf SOAR procedure}\label{soar1}
\begin{description}
  \item[1:]$q_1=u/\|u\|$, $p_1=0$                          \\[-8mm]
  \item[2:] for $j=1,2,\ldots,k$ do               \\[-8mm]
  \item[3:]\hspace{0.5cm} $r=Aq_j+Bp_j$           \\[-9mm]
  \item[4:]\hspace{0.5cm}  $s=q_j$                \\[-8mm]
  \item[5:]\hspace{0.5cm} for $i=1, 2, \ldots, j$ do \\[-9mm]
  \item[6:]\hspace{1cm}  $t_{ij}=q_i^*r$             \\[-9mm]
  \item[7:]\hspace{1cm}   $r=r-q_it_{ij}$            \\[-9mm]
  \item[8:]\hspace{1cm}  $s=s-p_it_{ij}$             \\[-8mm]
  \item[9:] \hspace{0.5cm} end for                  \\[-8mm]
  \item[10:] \hspace{0.5cm} $t_{j+1j}=\|r\|$         \\[-8mm]
  \item[11:]\hspace{0.5cm}   if $t_{j+1j}=0$, stop   \\[-8mm]
  \item[12:]\hspace{0.5cm }    $q_{j+1}=r/t_{j+1j}$  \\[-8mm]
  \item[13:]\hspace{0.5cm }    $p_{j+1}=s/t_{j+1j}$  \\[-8mm]
  \item[14:] end for                                 \\[-8mm]
\end{description}
\end{alg}

The following basic results hold for this algorithm; see \cite{bs}.

\begin{theorem}\label{thm1}
Define $Q_k=[q_1,q_2,\ldots,q_k]$ and $P_k=[p_1,p_2,\ldots,p_k]$
and $\hat{T}_k=\left[\begin{array}{c} T_k\\
t_{k+1k}e^*_k\end{array}\right]=[t_{ij}]\in {\cal C}^{(k+1)\times k}$.
If Algorithm~\ref{soar1} does not stop before step $k$, then we have
\begin{equation}
span\{Q_k\}=\mathcal{G}_k(A, B;u) \label{soarsub}
\end{equation}
and the $k$-step SOAR decomposition
\begin{equation}
H\left[
\begin{array}{c}Q_k \\P_k\end{array}\right]=
\left[ \begin{array}{c}Q_{k+1}
\\P_{k+1}\end{array}\right]\hat{T}_k,\label{05-6}
\end{equation}
where
$Q_{k+1}=[Q_k,q_{k+1}],~P_{k+1}=[P_k,p_{k+1}]$.
\end{theorem}

Before proceeding, we introduce the following definition.

\begin{defi}{\rm \cite{bs}}\label{def2}
If $r_i,\, i=0,1, \ldots, j$
are linearly dependent but $[r_i^T,r_{i-1}^T]^T,\,i=0, \ldots, j$
with $r_{-1}=0$ are not, we call this situation deflation; if
both $\{r_i\}$ and ~$\{[r_i^T,r_{i-1}^T]^T\}$ are
linearly dependent at step $j$, we call this situation breakdown.
\end{defi}

According to Definition~\ref{def2}, if Algorithm~\ref{soar1} stops
prematurely at step $j<k$, then either deflation or breakdown must
occur at that step. Deflation means that
$\mathcal{G}_{j+1}(A,B;u)=\mathcal{G}_j(A, B;u)$ but
${\cal K}_{j+1}(H,v)\not={\cal K}_j(H,v)$,
so the Arnoldi process on $H$ does not terminate at step $j$. As a result,
when deflation occurs at step $j$,
${\cal K}_j(H,v)$ does not contain any exact eigenvector of $H$,
which, from (\ref{gk1}), implies that ${\cal G}_j(A,B;u)$ may not
contain any exact eigenvector of QEP (\ref{05-qep}). Therefore,
deflation must be remedied to continue the algorithm.

Bai and Su \cite{bs} present the following algorithm that
detects and remedies deflation.

\begin{alg} {\bf SOAR procedure with deflation remedy} \label{soar2} \\[-8mm]
\begin{description}
  \item[1:]$q_1=u/\|u\|$, $p_1=0$                               \\[-8mm]
  \item[2:] for $j=1,2,\ldots,k$ do                      \\[-8mm]
  \item[3:]\hspace{0.5cm} $r=Aq_j+Bp_j$                  \\[-9mm]
  \item[4:]\hspace{0.5cm}  $s=q_j$                       \\[-8mm]
  \item[5:]\hspace{0.5cm} for $i=1, 2, \ldots, j$ do     \\[-9mm]
  \item[6:]\hspace{1cm}  $t_{ij}=q_i^*r$                 \\[-8mm]
  \item[7:]\hspace{1cm}   $r=r-q_it_{ij}$                \\[-9mm]
  \item[8:]\hspace{1cm}  $s=s-p_it_{ij}$                 \\[-9mm]
  \item[9:] \hspace{0.5cm} end for                      \\[-9mm]
  \item[10:] \hspace{0.5cm} $t_{j+1j}=\|r\|$             \\[-8mm]
  \item[11:]\hspace{0.5cm}   if $t_{j+1j}=0$             \\[-9mm]
  \item[12:]\hspace{1cm}    if $s\in span\{p_i|i:
                             q_i=0,1\leq i\leq j  \}$    \\[-8mm]
  \item[13:]\hspace{1.5cm}   break                       \\[-9mm]
  \item[14:]\hspace{1cm}   else   {\em deflation}        \\[-9mm]
  \item[15:]\hspace{1.5cm} reset $t_{j+1j}=1$            \\[-9mm]
  \item[16:]\hspace{1.5cm}  $q_{j+1}=0$                  \\[-8mm]
  \item[17:]\hspace{1.5cm}   $p_{j+1}=s$                 \\[-9mm]
  \item[18:]\hspace{1cm}     end if                      \\[-9mm]
  \item[19:]\hspace{0.5cm}   else                        \\[-9mm]
  \item[20:]\hspace{1cm }    $q_{j+1}=r/t_{j+1j}$        \\[-9mm]
  \item[21:]\hspace{1cm }    $p_{j+1}=s/t_{j+1j}$        \\[-9mm]
  \item[22:]\hspace{0.5cm }   end if                     \\[-9mm]
  \item[23:]   end for                                   \\[-8mm]
\end{description}
\end{alg}
In the procedure, if deflation
occurs, we simply set $t_{j+1 j}$ to one  and take $q_{j+1} = 0$.
To decide if $s\in span\{p_i|i: q_i=0,1\leq i\leq j \}$, the
Gram--Schmidt orthogonalization with refinement is used \cite{bs,oc}. When
deflation occurs, the nonzero vectors in the sequence $\{q_j\}$ are still
orthonormal and span the second-order Krylov subspace
$\mathcal{G}_k(A,B;u)$ with the dimension smaller than $k$.
We refer the reader to Bai and Su \cite{bs} for details.

We point out that Theorem~\ref{thm1} is true for Algorithm~\ref{soar2}
but there are zero columns in $Q_k$ when deflation occurs.



It is easily checked that a serious disadvantage of the SOAR procedure is that the
implicit restarting technique is not applicable since the updated $p_1$
is not zero any more. Several researchers have proposed replacing
$p_1=0$ in Algorithms~\ref{soar1}--\ref{soar2} by a nonzero one \cite{bao,oc,zhou}.
This leads to the following generalized second-order
Krylov sequence and subspace; see \cite{bao,zhou}.

\begin{defi}  \label{def3}
Let $A$ and $B$ be $n\times n$ matrices and for vectors $u_1, u_2\in
\mathcal{C}^n$, and define
$$
\begin{array}{ccl}
r_0&=&u_1,\\
r_1&=&Ar_0+Bu_2,\\
r_j&=&Ar_{j-1}+Br_{j-2} \hspace{0.5cm} \mbox{for } j\geq 2.
\end{array}
$$
Then $r_0, r_1, r_2, \ldots, r_{k-1}$
is called a generalized second-order Krylov sequence
based on $A,B$ and $u_1,u_2$,
and
$
\mathcal{G}_k(A,B;u_1,u_2)={span}\{r_0, r_1, r_2, \ldots, r_{k-1}\}
$
the $k$-th generalized second-order Krylov subspace.
\end{defi}

Obviously, $\mathcal{G}_k(A,B;u_1,0)=\mathcal{G}_k(A,B;u_1)$.
For a general $\tilde v=[u_1^T, u_2^T]^T$, it is seen that
\begin{equation}
\left[\begin{array}{c}r_j\\r_{j-1} \end{array} \right]=H^j\tilde v, \ j\geq 1.
\label{newlinsub}
\end{equation}
But different from (\ref{gk1}), since $u_2$ is a general vector,
the fundamental relation now becomes
\begin{equation}
H{\cal K}_{k-1}(H,\tilde v)=span\{H\tilde v,\ldots,H^{k-1}\tilde v\}
\subseteq{\cal G}_k^2(A,B;u_1,u_2), \label{gk2}
\end{equation}
where ${\cal G}_k^2(A,B;u_1,u_2)$ is the subspace generated by the vector set
$$
\left\{\left[\begin{array}{c}
r_0\\
0\\
\end{array}
\right],
\left[\begin{array}{c}
r_1\\
0\\
\end{array}
\right]
,\ldots,
\left[\begin{array}{c}
r_{k-1}\\
0\\
\end{array}
\right],
\left[\begin{array}{c}
0\\
r_0\\
\end{array}
\right],
\left[\begin{array}{c}
0\\
r_1\\
\end{array}
\right],
\ldots,
\left[\begin{array}{c}
0\\
r_{k-1}\\
\end{array}
\right]
\right\}.
$$
Note that if the eigenvector $[\lambda x^T,x^T]^T$ is
contained in ${\cal K}_{k-1}(H,\tilde v)$ then it also lies in the subspace
$H{\cal K}_{k-1}(H,\tilde v)$. If this is the case, (\ref{gk2}) shows
that the eigenvector $x$ of QEP~(\ref{05-qep}) is contained in
${\cal G}_k(A,B;u_1,u_2)$. More generally, by continuity, it is deduced from
the above that if ${\cal K}_{k-1}(H,\tilde v)$ has a good approximation to
$[\lambda x^T,x^T]^T$ then $H{\cal K}_{k-1}(H,\tilde v)$ has one too,
which, in turn, means that there must be a good
approximation to $x$ contained in ${\cal G}_k(A,B;u_1,u_2)$.

Analogous to Algorithm~\ref{soar2}, we can present a GSOAR procedure,
i.e., Algorithm~\ref{GSOAR}, that remedies deflation and generates the vector
sequence $\{q_j\}$, whose nonzero ones form an orthonormal basis of
$\mathcal{G}_k(A,B;u_1,u_2)$. We point out that
the GSOAR procedure in \cite{bao,zhou} is the same as Algorithm~\ref{soar1}
except that $p_1=0$ in line 1 is replaced by a general vector $p_1=u_2/\|u_2\|$.

\begin{alg} {\bf GSOAR procedure with deflation remedy}\label{GSOAR}
\begin{description}
  \item[1:] $q_1=\frac{u_1}{\|u_1\|}$, $p_1=\frac{u_2}{\|u_2\|}$.  \\[-7mm]
  \item[2:]for $j=1,2,\ldots,k$ do                                \\[-8mm]
  \item[3:] \hspace{0.5cm} $r=Aq_j+Bp_j$                 \\[-8mm]
  \item[4:]\hspace{0.5cm}  $s=q_j $                      \\[-7mm]
  \item[5:]\hspace{0.5cm}  for ~$i=1,2,\ldots,j$ do               \\[-8mm]
  \item[6:]\hspace{1cm}     $t_{ij}=q_i^*r $             \\[-7mm]
  \item[7:]\hspace{1cm}     $r=r-t_{ij}q_i$              \\[-8mm]
  \item[8:]\hspace{1cm}     $s=s-t_{ij}p_i$              \\[-7mm]
  \item[9:]\hspace{0.5cm}  end for                       \\[-8mm]
  \item[10:]\hspace{0.5cm}   $t_{j+1j}=\|r\|$             \\[-7mm]
  \item[11:]\hspace{0.5cm}  if  ~$t_{j+1j}=0$             \\[-8mm]
  \item[12:]\hspace{1cm}    if $s\in span\{p_i|i:
                             q_i=0,1\leq i\leq j  \}$    \\[-8mm]
  \item[13:]\hspace{1.5cm}   break                       \\[-9mm]
  \item[14:]\hspace{1cm}   else   {\em deflation}        \\[-9mm]
  \item[15:]\hspace{1.5cm} reset $t_{j+1j}=1$            \\[-9mm]
  \item[16:]\hspace{1.5cm}  $q_{j+1}=0$                  \\[-8mm]
  \item[17:]\hspace{1.5cm}   $p_{j+1}=s$                 \\[-9mm]
  \item[18:]\hspace{1cm}     end if                      \\[-9mm]
  \item[19:]\hspace{0.5cm}   else                        \\[-9mm]
  \item[20:]\hspace{1cm }    $q_{j+1}=r/t_{j+1j}$        \\[-9mm]
  \item[21:]\hspace{1cm }    $p_{j+1}=s/t_{j+1j}$        \\[-9mm]
  \item[22:]\hspace{0.5cm }   end if                     \\[-9mm]
  \item[23:]   end for                                   \\[-8mm]
\end{description}
\end{alg}

It is direct to justify that Theorem~\ref{thm1} holds for this algorithm
with a general $p_1=u_2/\|u_2\|$, that is, we have
$$
span\{Q_k\}={\cal G}_k(A,B;q_1,p_1)
$$
with $q_1$ and $p_1$ normalized
and the $k$-step GSOAR decomposition (\ref{05-6}) if the algorithm
does not break down before step $k$.

Otto \cite{oc} defines the modified second-order Krylov sequence as
$r_0, r_1, r_2, \ldots, r_{k-1}, u_2$ and the modified second-order Krylov
subspace of dimension $k+1$ generated by the vector sequence.
After the orthonormal $q_1,q_2,\ldots,q_{k+1}$ are
generated, he orthonormalizes $u_2$ against them to
get $q_{k+2}$. This is called the modified SOAR procedure. A disadvantage of
it is that there is no compact relationship (\ref{gk1})
or (\ref{gk2}). So it is hard to interpret such a modified subspace
and establish definitive results on breakdown and deflation.

The GSOAR method is a Rayleigh--Ritz method, and it projects the large
QEP (\ref{05-qep}) onto ${\cal G}_k(A,B;u_1,u_2)$ by imposing
the Galerkin condition, leading to the $k$-dimensional QEP
\begin{equation}
(\theta^2M_k+\theta C_k+K_k)g=0 \label{05-pro}
\end{equation}
with $\|g\|=1$, where  $M_k=Q_k^*MQ_k$, $C_k=Q_k^*CQ_k$ and $K_k=Q_k^*KQ_k$.
Let the $(\theta, g)$ be the eigenpairs of (\ref{05-pro}). Then
the GSOAR method uses the Rayleigh--Ritz pairs $(\theta, y(=Q_kg))$
to approximate some of the eigenpairs of (\ref{05-qep}). We comment that if
deflation occurs then $Q_k$ consists of only nonzero orthonormal vectors
$q_j$ and the dimension of (\ref{05-pro}) is smaller than $k$.

\section{A refined GSOAR (RGSOAR) method}\label{s3}

As is known, the Rayleigh--Ritz method may fail to converge for computing
eigenvectors of the linear eigenvalue problem and the QEP; see \cite{jiastewart}
and \cite{huangjialin}, respectively. To correct this deficiency, a refined
projection principle is proposed in \cite{J44} (see also \cite{stewart,vorst})
for the linear eigenvalue problem, which leads to the refined Rayleigh--Ritz method.
The refined method extracts the best approximate eigenvectors
from a given subspace in the sense that the residuals formed with certain approximate
eigenvalues available are minimized in the sense of 2-norm over the subspace.
A refined GSOAR (RGSOAR) method has been proposed in \cite{bao,zhou}. We next describe it
and give more details on some practical issues.

Suppose that we have computed the Ritz values $\theta$ by the
GSOAR method and select $m$ ones of them to approximate
$m$ desired eigenvalues of (\ref{05-qep}). For each chosen $\theta$,
the RGSOAR method seeks a unit length vector
$\tilde{u} \in \mathcal{G}_k(A, B;u_1,u_2)$
satisfying the optimal requirement
\begin{equation}
\|(\theta^2 M + \theta C+ K) {\tilde{u}}\| = \min_{
 \mbox{\scriptsize $\begin{array}{c} u \in \mathcal{G}_k(A,
B;u_1,u_2) \\ \|u\|=1
 \end{array}$}}
\|(\theta^2 M + \theta C+ K )u\| \label{th2}
\end{equation}
and uses it as an approximate eigenvector, called the refined Ritz
vector. The pairs $(\theta,\tilde u)$ are also called the refined Rayleigh--Ritz
approximations. Since the (non-zero) columns of $Q_k$ form an orthonormal
basis of $\mathcal{G}_k(A, B;u_1,u_2)$,
(\ref{th2}) amounts to seeking a unit length vector
$\tilde{z} \in {\cal C}^k$ such that
$\tilde{u}=Q_k\tilde z$ with
\begin{equation}
\tilde z=\arg\min_{ \mbox{\scriptsize $\begin{array}{c} z \in {\cal C}^k \\
\|z\|=1 \end{array}$}} \|(\theta^2 M + \theta C+ K)Q_kz\|,
\label{th3}
\end{equation}
the right singular vector of the matrix
$\theta^2 MQ_k+\theta CQ_k+KQ_k$ associated with its smallest
singular value $\sigma_{\min}
(\theta^2 MQ_k+\theta CQ_k+KQ_k)$. However, the direct computation
of its SVD may be expensive. Precisely, assume that the matrix is real
and $k\ll n$. Then the cost of Golub--Reinsch's SVD
algorithm is about $4nk^2$ flops, and that of Chan's SVD
algorithm is about $2nk^2$ flops \cite[p. 254]{golub}. Keep in mind
that $m$ is the number of the desired eigenpairs. The CPU time costs are
then $4nmk^2$ and $2nmk^2$ flops, respectively.

The first author in \cite{jsvd} has proposed a cross-product matrix-based
algorithm for computing the SVD of a matrix, which can be much more
efficient than the above standard SVD algorithms. Applying the
algorithm to (\ref{th2}), we form the cross-product matrix
$$B_k=\left(\theta^2 MQ_k+\theta CQ_k+KQ_k \right )^*
\left(\theta^2 MQ_k+\theta CQ_k+KQ_k \right),
$$
which is the Hermitian (semi-)positive definite. $\tilde z$ is then the
eigenvector of $B_k$ associated with its smallest eigenvalue $\sigma^2_{\min}
(\theta^2 MQ_k+\theta CQ_k+KQ_k)$. We compute the eigensystem
of $B_k$ by the QR algorithm to get $\tilde z$. In finite precision arithmetic,
the computed eigenvector is an approximation to $\tilde{z}$ with accuracy
$O(\epsilon_{\rm mach})$ provided that the second smallest
singular value of $\theta^2 MQ_k+\theta CQ_k+KQ_k$ is not very close to
the smallest one, where $\epsilon_{\rm mach}$ is the machine precision.

Let us now look at the computational cost of this algorithm. Define
$$
W_1=MQ_k,\,W_2=CQ_k,\,W_3=KQ_k,
$$
which are available when forming the projected QEP and do not
need extra cost. Then
\begin{eqnarray}
B_k&=&\mid\theta\mid^4W_1^*W_1+\mid\theta\mid^2W_2^*W_2+W_3^*W_3 +
\theta\bar{\theta}^2W_1^*W_2+\bar{\theta}\theta^2W_2^*W_1\\\nonumber
&&+\bar{\theta}^2W_1^*W_3+\theta^2W_3^*W_1+\bar{\theta}W_2^*W_3+\theta
W_3^*W_2,
\end{eqnarray}
where the bar denotes the complex conjugate of a scalar. Assume that
$W_1,W_2$ and $W_3$ are real and note that $B_k$ is Hermitian for a
complex $\theta$ and real symmetric for a real $\theta$. Then we only need to
form the upper (lower) triangular part of $B_k$, which involves the
upper (lower) triangular parts of the nine matrices
$W_i^*W_j,\,i,j=1,2,3$. All these cost about $9nk^2$ flops. With
these nine $W_i^*W_j$ available, we only need $O(k^2)$ flops to form $B_k$
for either a real or complex $\theta$, negligible to $9nk^2$ flops.
So, we CPU timely need $9nk^2$ flops to form $m$ Hermitian matrices $B_k$ for $m$
approximate eigenvalues $\theta$. We then compute the complete
eigensystems of these $B_k$ by the QR algorithm using $O(mk^3)$
flops. Therefore, we can compute $m$ right singular vectors
$\tilde z$ using about $9nk^2$ flops when $mk\ll n$, a natural requirement
in practice. As a result, a simple comparison indicates that such cross-product
based algorithm is more efficient than Golub--Reinsch's SVD algorithm when
$m\geq 3$ and Chan's SVD algorithm when $m\geq 5$.

We can now present a basic (non-restarted) RGSOAR algorithm.

\begin{alg} {\bf The RGSOAR algorithm}
\label{RGSOAR}
\begin{enumerate}
\item Given the starting vectors $u_1,u_2$,
run the GSOAR procedure to generate an
orthonormal basis $Q_k$ of
$\mathcal{G}_k(A,B;u_1,u_2)$.

\item Compute $W_1=MQ_k,~W_2=CQ_k$ and $W_3=KQ_k$.

\item Compute $M_k=Q_k^*W_1,~C_k=Q_k^*W_2$ and $K_k=Q_k^*W_3$,
solve the projected QEP
  \begin{equation} (\theta_i^2M_k+\theta_iC_k+K_k)g_i=0,
  \label{refined-QEP}
  \end{equation}
and select $m$ Ritz values $\theta_i$ as
approximations to the $m$ desired eigenvalues $\lambda_i$.

\item For each chosen $\theta_i,\,1\leq i\leq m$, form $B_k$, and compute the
  eigenvector $\tilde{z}_i$ of $B_k$ associated with its smallest
  eigenvalue and the refined Ritz vector $\tilde{u}_i=Q_k\tilde{z}_i$.

\item Test convergence of $(\theta_i,\tilde{u}_i)$ by
computing the relative residual norms
$$
\frac{\|(\theta_i^2M+\theta_i C+K)\tilde{u}_i\|}{|\theta_i|^2\|M\|_1
+|\theta_i|\|C\|_1+\|K\|_1},\ i=1,2,\ldots,m.
$$
\end{enumerate}
\end{alg}

\section{Implicitly restarted algorithms}\label{s4}

This section consists of three subsections.  In Section~\ref{s4-2},
under the assumption that no deflation occurs, we describe how
to implicitly restart the GSOAR procedure. In Section~\ref{s4-3}, we discuss
how to select best possible shifts, and propose exact and
refined shifts for respective use within implicitly restarted GSOAR and RGSOAR
algorithms. In Section~\ref{s4-4}, we present an effective approach to cure deflation
in implicit restarts, so that implicit restarting can be run unconditionally.

\subsection{Implicit restarts}\label{s4-2}

As step $k$ increases, the GSOAR and RGSOAR methods become expensive and
impractical due to storage requirement and/or computational
cost. So restarting is generally necessary. That is, for a given maximum
$k$ and the subspace
${\cal G}_k(A,B;q_1,p_1)$ with $q_1$ and $p_1$ normalized, if the methods do
not converge yet, based on the information available, we select new unit length
vectors $q_1^+$ and $ p_1^+$ to construct a better subspace
$\mathcal{G}_k(A, B;q_1^+,p_1^+)$ that contains richer information on the
desired eigenvectors $x$. We then extract new better
approximate eigenpairs with respect to $\mathcal{G}_k(A, B;q_1^+,p_1^+)$.
Proceed in such a way until the methods converge.

If no deflation occurs, it is direct to adapt the implicit restarting scheme
\cite{so} to the modified SOAR procedure in \cite{oc} and the GSOAR
procedure. Given $p$ shifts $\mu_1,\mu_2,\ldots,\mu_p$,
performing $p$ implicit shifted QR iterations on $T_k$ yields the relation
$$
(T_k-\mu_1I)\cdots (T_k-\mu_pI)=V_k R,
$$
where $V_k$ is a $k\times k$ orthogonal (unitary) matrix and $R$ is
upper triangular. Specifically, $V_k$ has only $p$ nonzero
subdiagonals. Adapted from the derivation of implicitly
restarting the standard Arnoldi process \cite{so}, we can establish
the following result for the GSOAR procedure.

\begin{theorem}\label{restart}
Given $p$ shifts $\mu_1,\ldots,\mu_p$, perform $p$
steps of implicit shifted QR iterations on $T_k$.
Let $\psi(T_k)=V_kR_k$ with $\psi(\mu)=\prod_{j=1}^p(\mu-\mu_j)$, and
define $Q_k^+=Q_kV_k$ and
$T_k^+=V_k^*T_kV_k$. Assume that no deflation occurs in the
$k(=m+p)$-step GSOAR decomposition {\rm (\ref{05-6})}. Then we
have an updated $m$-step GSOAR decomposition
\begin{equation}
H \left[
\begin{array}{c}Q_m^{+} \\P_m^{+}\end{array}\right]=
\left[ \begin{array}{c}Q_m^{+}
\\P_m^{+}\end{array}\right]T_m^{+}+
\tilde{t}_{m+1m}^{+}\left [ \begin{array}{c}
q_{m+1}^{+}\\p_{m+1}^{+}\end{array} \right] e^*_m \label{resrt-3}
\end{equation}
starting with $\left[
\begin{array}{c}q_1^+ \\p_1^+\end{array}\right]$,
where $Q_m^+=Q_kV_k(:,1:m)$, $P_m^+=P_kV_k(:,1:m)$,
$T_m^+=T_k^+(1:m,1:m)$ is upper Hessenberg and
\begin{eqnarray*}
\left [
\begin{array}{c} q_{m+1}^{+}\\p_{m+1}^{+}\end{array} \right]&=&
\frac{1}{\tilde{t}_{m+1m}^{+}}f_m^+,\\
f_m^+&=&{t}_{m+1m}^+\left[ \begin{array}{c}q_{m+1}^+ \\
p_{m+1}^+\end{array}\right]+t_{k+1k}V_k(k,m)\left [
\begin{array}{c} q_{k+1}\\p_{k+1}\end{array} \right],
\\
\tilde{t}_{m+1m}^+&=&\|{t}_{m+1m}^+ q_{m+1}^+ +t_{k+1k}V_k(k,m)q_{k+1}\|
\end{eqnarray*}
with $V_k(k,m)$ the entry of $V_k$ in position $(k,m)$.
\end{theorem}

Theorem~\ref{restart} states that if no deflation occurs then we have
naturally obtained an $m$-step GSOAR decomposition (\ref{resrt-3}) after $p$
implicit shifted QR iterations are run on $T_k$, thus generating an orthonormal
basis $\{q_j\}_{j=1}^m$ of the $m$-dimensional subspace ${\cal G}_m(A,B;
q_1^+,p_1^+)$. Decomposition (\ref{resrt-3}) is then extended to a $k$-step one
from step $m+1$ upwards in a standard way other than from scratch, producing
an orthonormal basis $\{q_j^+\}_{j=1}^k$ of the updated $k$-dimensional
subspace ${\cal G}_k(A,B;q_1^+,p_1^+)$.

Analogous to the proof of the result on updated starting vectors
in \cite{so}, it is direct to justify the following theorem.

\begin{theorem}\label{revector}
It holds that
\begin{equation}\label{updatevc}
\left[
\begin{array}{c}q_1^+ \\p_1^+\end{array}\right]
=\frac{1}{\tau }\psi(H)\left[
\begin{array}{c}q_1 \\p_1\end{array}\right],
\end{equation}
with $\psi(\lambda)=\prod_{j=1}^p(\lambda - \mu_j)$
and $\tau$ a normalizing factor.
\end{theorem}

\subsection{The selection of shifts}\label{s4-3}

The selection of the shifts is one of the keys for the success and overall
efficiency of implicitly restarted GSOAR and RGSOAR algorithms.
In this subsection we propose the corresponding best possible shifts for
respective use within each algorithm.

Assume that $H$ is diagonalizable.
It is shown in, e.g., \cite{saad1}, that if the starting vector $\tilde v$
is a linear combination of $m$ eigenvectors of $H$ then ${\cal K}_m(H,\tilde v)$
is an invariant subspace. Therefore, a fundamental principle
of restarting is to select a better
vector $\tilde v^+$, in some sense, from the current ${\cal K}_k(H,\tilde v)$ as an
updated starting vector that amplifies the components of the desired eigenvectors
and simultaneously dampens those of the unwanted ones,
so that the updated ${\cal K}_k(H,\tilde v^+)$ contains
more accurate approximations to the $m$ desired eigenvectors.
For implicit restarting, based on formulas for updated
starting vectors like (\ref{updatevc}), for the linear eigenvalue
problem and the computation of a partial SVD, it has been shown in \cite{J50,jia02}
and \cite{jianiu03,jianiu10} that such goal is achieved by selecting the shifts to
approximate some of the unwanted eigenvalues or singular values as best as possible
within the framework of the underlying method.
A general result is that the better the shifts approximate the
unwanted eigenvalues, the richer information on the desired eigenvectors is
contained in the updated starting vector, so that a better Krylov
subspace is generated.

Motivated by the above results, we now investigate a reasonable selection of
shifts for use within implicitly restarted GSOAR and RGSOAR algorithms.
Observe that the projected QEP (\ref{05-pro}) of the
large QEP (\ref{05-qep}) over $span\{Q_k\}$ amounts to the generalized
eigenvalue problem
\begin{equation}\label{projgep}
\left[\begin{array}{cc}
-C_k &-K_k\\
I & 0
\end{array}\right]
\left[\begin{array}{c}
\theta g\\
g\end{array}\right]=\theta \left[\begin{array}{cc}
M_k & 0\\
0 & I
\end{array}\right]\left [\begin{array}{c}
\theta g\\
g\end{array}\right ],
\end{equation}
which is the projected problem of large generalized eigenvalue problem
(\ref{5xian}) over the subspace ${\cal G}^2_k(A,B;u_1,u_2)$ (c.f. (\ref{gk2}))
spanned by the (nonzero) columns of
$$
\hat{Q}_{2k}=\left[\begin{array}{cc}
Q_k & 0\\
0 &Q_k
\end{array}
\right].
$$
The above problem amounts to the standard linear eigenvalue problem
$$
\left[\begin{array}{cc}
-M_k^{-1}C_k &-M_k^{-1}K_k\\
I & 0
\end{array}\right]
\left[\begin{array}{c}
\theta g\\
g\end{array}\right]=\theta \left[\begin{array}{c}
\theta g\\
g\end{array}\right].
$$
(\ref{updatevc}) indicates that we should select
the shifts $\mu_j,\ j=1,2,\ldots,m$ as the best possible approximations to the
unwanted eigenvalues of $H$ so as to generate increasingly better updated
subspaces ${\cal K}_k(H,\tilde v^+)$ and $H{\cal K}_{k-1}(H,\tilde v)$  with
$\tilde v^+=[{q_1^+}^T,{p_1^+}^T]^T$. In terms of (\ref{gk2}) and the comments
followed, this, in turn, leads to increasingly better updated
${\cal G}_k^2(A,B;q_1^+,p_1^+)$  that contains increasingly better
approximations to the $m$ desired eigenvectors of $H$. As a result,
${\cal G}_k(A,B;q_1^+,p_1^+)$ contains more accurate approximations to
the desired eigenvectors of (\ref{05-qep}).
So, just as for the linear eigenvalue problem, we
should choose shifts for each implicitly restarted GSOAR type algorithm
in the sense that they are best possible approximations to some of the unwanted
eigenvalues of (\ref{05-qep}).

For the Rayleigh--Ritz method with respect to a given subspace,
the Ritz values can be considered as the best
approximations available to some eigenvalues of (\ref{05-qep}).
Otto \cite{oc} proposed exact second-order shifts for his implicitly restarted
modified SOAR algorithm. Adapted here, one solves the projected QEP (\ref{05-pro})
and selects $m$ Ritz values $\theta_i$ as approximations to the desired eigenvalues.
Then the unwanted Ritz values are shift candidates, called the
exact second-order shift candidates. A problem is that there are $2k-m$
shift candidates, while for (\ref{resrt-3}) the
number $p$ of shifts must not exceed $k-m$. One must select $p=k-m$ shifts among
the $2k-m$ candidates. Otto simply suggested to take any $p=k-m$ shifts
among $2k-m$ ones. We should point out that this situation is unlike implicitly
restarted Arnoldi type algorithms for the linear eigenvalue problem, where
the the maximum number of shifts is just that of candidates; see \cite{so}
and \cite{J50,jia02,jianiu03,jianiu10}.

However, the above selection of exact second-order shifts is
problematic and susceptible to failure, as elaborated below.
It is crucial to keep in mind a basic fact that the QEP may often have two
distinct eigenvalues that share the same eigenvector \cite{tm}.
This means that, for QEP (\ref{05-pro}),
some of the shift candidates and some of the $m$ Ritz values used to approximate
the desired eigenvalues may share common eigenvector(s). Therefore, if it is
unfortunate to take such candidates for shifts, restarting will filter out the
information on the corresponding desired eigenvectors
and thus makes implicitly restarted GSOAR algorithms perform poorly.

In order to avoid the above deficiency, we propose new shift candidates
for the implicitly restarted GSOAR and RGSOAR
algorithms, respectively, and show how to reasonably
select the shifts among the candidates. We first consider the GSOAR method.
Project QEP (\ref{05-qep}) onto the orthogonal complement of
$span\{{y}_1,\ldots,{y}_m\}$ with respect to
$\mathcal{G}_k(A,B;q_1,p_1)$, where ${y}_1,\ldots,{y}_m$ are the Ritz
vectors approximating the desired eigenvectors $x_1,\ldots,x_m$.
Then we obtain a $p$-dimensional projected QEP and compute its
$2p$ eigenvalues. A remarkable consequence is that these $2p$ eigenvalues must be
approximations to some of the unwanted eigenvalues of QEP
(\ref{05-qep}) because the information on $x_1,\ldots,x_m$ has been
removed from $\mathcal{G}_k(A,B;q_1,p_1)$.  So we can
use any $p$ ones of these $2p$ candidates as shifts. To be unique, we choose the $p$
ones {\em farthest} from the Ritz values $\theta_i,\,i=1,2,\ldots,m$ that are
used to approximate the desired eigenvalues $\lambda_1,\ldots,\lambda_m$.
The motivation of this choice is that, based on (\ref{updatevc}),
these shifts can be better to amplify the information of $\tilde v^+$ on
the desired eigenvectors and dampen the components of undesired eigenvectors
in $\tilde v^+$.

If we are interested in the $m$ eigenvalues nearest to a target $\sigma$
and/or the associated eigenvectors, QEP~(\ref{05-qep}) can be
equivalently transformed to a shift-invert QEP; see the end of
this subsection. In this case, we select
the $p$ Ritz values among $2p$ candidates {\em farthest} from $\sigma$ as shifts.
Such selection of shifts is motivated by an idea from \cite{jianiu03,jianiu10}, where
some of the shifts are taken to be unwanted Ritz values farthest from the wanted
approximate singular values. It was argued there that this selection
can better dampen those components of the unwanted
singular vectors and meanwhile amplify the components of
the desired singular vectors.

We now turn to the selection of shifts for the RGSOAR algorithm.
Algorithm \ref{RGSOAR} computes the refined Ritz vectors
$\tilde{u}_i$, which are generally more and can be much more accurate
than the Ritz vectors $y_i$ \cite{huangjialin,jiastewart}.
The first author \cite{J50,jia02} has proposed certain refined shifts
for the refined Arnoldi method and the refined harmonic Arnoldi method for the
linear eigenvalue problem. It is shown that the refined shifts are
generally better than the corresponding exact shifts and can be
computed efficiently and reliably. In the same spirit, we next propose certain
refined shifts for the RGSOAR algorithm.

Since the refined Ritz vectors $\tilde{u}_i,\,i=1,2,\ldots,m$ are more accurate
than the corresponding $y_i$, the orthogonal complement of
$span\{\tilde{u}_1,\ldots,\tilde{u}_m\}$ with respect to
$\mathcal{G}_k(A,B;q_1,p_1)$ contains richer information on the
unwanted eigenvectors than the orthogonal
complement of $span\{y_1,\ldots,y_m\}$ with respect to
$\mathcal{G}_k(A,B;q_1,p_1)$. As a result, the
eigenvalues of the projected QEP of QEP (\ref{05-qep}) onto this
orthogonal complement are more accurate approximate eigenvalues than
the exact shift candidates described above. We call them refined shift candidates.
We use the same approach as above to select $p$ ones among them as shifts,
called the refined shifts, for use within the implicitly restarted RGSOAR algorithm.

Finally, we show how to compute the exact and refined shifts efficiently
and reliably. We take the refined shifts as example. The computation
of exact shifts is analogous. Recall $\tilde{u}_i =Q_k \tilde{z}_i,
i=1,2,\ldots ,m$, and write $Z_m=[\tilde{z}_1,\ldots ,\tilde{z}_m]$.
If QEP (\ref{05-qep}) is real and two columns $\tilde{z}_i$ and $\tilde{z}_{i+1}$
of $Z_m$ are complex conjugate, we replace them by their normalized real
and imaginary parts, respectively, so that the resulting $Z_m$ is real.
We then make the full QR decomposition
$$Z_m=[U_m ,U_{\perp}]
 \left[ \begin{array}{c} R_m  \\
      0   \end{array} \right],$$
where $U_m$ and $U_{\perp}$ are $k \times m$ and $k\times p$ column
orthonormal matrices, respectively, and $R_m$ is $m\times m$
upper triangular. We use the Matlab built-in function
{\sf qr.m} to compute the decomposition in experiments. This costs $O(k^3)$ flops,
negligible to the cost of the $k$-step GSOAR procedure.
Obviously, it holds that
$$
span\{\tilde u_1,\ldots,\tilde u_m\}=span\{Q_kU_m\},\ \ \
span\{[Q_kU_m,Q_kU_{\perp}]\}={\cal G}_k(A,B;q_1,p_1).
$$
Therefore,
$Q_kU_{\perp}$ is an orthonormal basis of the orthogonal
complement of $span\{\tilde{u}_1,\ldots,\tilde{u}_m\}$ with respect to
$\mathcal{G}_k(A,B;q_1,p_1)$. It is direct to justify that the projected QEP of the
original QEP (\ref{05-qep}) onto $span\{Q_kU_{\perp}\}$ is just the projected QEP of
the small QEP (\ref{refined-QEP}) onto $span\{U_{\perp}\}$. So, we
form the projected QEP of the original QEP (\ref{05-qep})
onto $span\{Q_kU_{\perp}\}$ at cost of
$O(k^3)$ flops. We then compute its $2p$ eigenvalues using $O(p^3)$
flops and select $p$ ones among them as the refined shifts. Since $p<k$,
the CPU time cost of computing the refined shifts is $O(k^3)$ flops.
For the exact shifts, recall the Ritz vectors
$y_i=Q_kg_i,\,i=1,2,\ldots,m$. Write $G_m=[g_1,\ldots,g_m]$ and
replace $Z_m$ by it. We then compute the exact shifts in the same way as above.

Having done the above, we have finally developed the following Algorithm~\ref{pracim}.

\begin{alg}{\bf The implicitly restarted GSOAR type
algorithms}\label{pracim}
\begin{enumerate}
\item Given unit length starting vectors $q_1$ and $p_1$, the number $m$
  of desired eigenpairs and the number $p$ of shifts $p$ satisfying $p\leq m-k$,
  run the $k$-step GSOAR procedure to generate $Q_k$.

\item Do until convergence

  Project QEP
            {\rm (\ref{05-qep})} onto $span\{Q_k\}$ to get QEP {\rm (\ref{05-pro})},
      select $m$ Ritz pairs $(\theta_i,y_i)$ or
            refined Ritz pairs $(\theta_i,\tilde{u}_i)$ as approximations to
            the $m$ desired eigenpairs, respectively, and determine their convergence.

\item If not converged, compute the $p$ exact shifts or refined shifts,
  and implicitly restart the GSOAR method or the
  RGSOAR method, respectively.

\item EndDo
\end{enumerate}
\end{alg}

Algorithm~\ref{pracim} includes two algorithms: the implicitly
restarted GSOAR algorithm with the exact shifts and
RGSOAR algorithm with the refined shifts, abbreviated as IGSOAR and
IRGSOAR here and hereafter. They can be used to compute a number of
largest eigenvalues in magnitude and the associated eigenvectors of
QEP (\ref{05-qep}). We determine the convergence of a Ritz pair
$(\theta,y)$ by requiring
\begin{equation}\label{stop}
\frac{\|(\theta^2M+\theta C+K)y\|}{|\theta|^2\|M\|_1+|\theta|\|C\|_1+
\|K\|_1}\leq tol,
\end{equation}
where $tol$ is a user-prescribed accuracy.
For the convergence of a refined Ritz pair $(\theta,\tilde{u})$, we
replace the above $y$ by $\tilde{u}$.

If the $m$ eigenvalues closest to a given target $\sigma$ are desired, we
use the shift-invert transformation $\rho=\frac{1}{\lambda-\sigma}$ with
$\det(Q(\sigma))\not=0$ to transform QEP~(\ref{05-qep}) to the new QEP
\begin{equation}\label{shiftqep}
Q_{\sigma}(\rho)x=(\rho^2M_{\sigma}+\rho C_{\sigma}+K_{\sigma})x=0,
\end{equation}
where $M_{\sigma}=\sigma^2M+\sigma C+K$ is nonsingular as $\det(M_{\sigma})=
\det(Q(\sigma))\not=0$, $C_{\sigma}=C+2\sigma M$,
$K_{\sigma}=M $. We then apply the previous analysis and algorithms
to (\ref{shiftqep}). Let $(\tilde{\rho},y)$ be an
approximate eigenpair (either a Ritz or refined Ritz pair) of
$Q_{\sigma}(\rho)x=0$ and $\hat{r}=Q_{\sigma}(\tilde{\rho})y$. Then
$(\frac{1}{\tilde\rho}+\sigma,y)$ is the corresponding approximate
eigenpair of $Q(\lambda)x=(\lambda^2M+\lambda C+K)x=0$. Define
$\tilde r=Q(\frac{1}{\tilde\rho}+\sigma)y$.
Then we obtain
\begin{eqnarray}
{\hat{r}}/{\tilde{\rho}^2}& =& (M_{\sigma} + C_{\sigma}/\tilde{\rho} +
K_{\sigma}/\tilde{\rho}^2)y\nonumber\\
              & =&(\sigma^2M+\sigma C+K +{(C+2\sigma M)}/{\tilde{\rho}}
                  +M/{\tilde{\rho}^2})y\nonumber\\
              & =& ((\frac{1}{\tilde{\rho}}+\sigma)^2M
              +(\frac{1}{\tilde{\rho}}+\sigma)C+K)y=Q(\frac{1}
              {\tilde{\rho}}+\sigma)y=\tilde{r},
\end{eqnarray}
from which it is direct to get the desired $\|\tilde{r}\|$ from $\|\hat{r}\|$
without computing $\tilde{r}$ explicitly.

We make a final note on Algorithm~\ref{pracim}. In previous discussions
and analysis, we have supposed $p=k-m$ previously. This is not mandatory.
In order to compute $m$ desired eigenpairs of (\ref{05-qep}),
the only restriction to $p$ is that $p\leq k-m$.
So the choice of $p$ is flexible and takes the form $p=k-(m+l)$ with $l$ a very small
nonnegative integer, as done in \cite{so} and \cite{J50,jia02,jianiu03,jianiu10},
where $l=3$ is often used.
We remark that different $p$ may have considerable effects on the overall
performance of the algorithms, but its choice can only be empirical.

\subsection{Cure of deflations in implicit restarts}\label{s4-4}

Theorem~\ref{restart} requires that no deflation occurs
in implicit restarts. If deflations occur at steps $m_1,m_2,\ldots,
m_j\leq k$, then the corresponding $j$ columns $q_{m_j}$ of $Q_{k}$ are
zeros. Denote by $\hat{Q}_k$ and $\hat{V}_k$ the matrices by
deleting the zero columns of $Q_k$ and rows $m_1,m_2, \ldots, m_j$
of $V_k$, respectively. Then we have $Q_{k}^{+}=Q_kV_k =
\hat{Q}_k\hat{V}_k$, from which and (\ref{05-6}) we get
\begin{equation}
\left[ \begin{array}{cc}A&B \\I & 0\end{array}\right] \left[
\begin{array}{c}\hat{Q}_k\hat{V}_k \\P_{k}V_k\end{array}\right]=
\left[ \begin{array}{c}\hat{Q}_k\hat{V}_k
\\P_{k}V_k\end{array}\right]T_{k}^+ +
t_{k+1k}\left [ \begin{array}{c} q_{k+1}\\p_{k+1}\end{array} \right]
e^*_{k}V_k, \label{defl-1}
\end{equation}
where $T_{k}^+=V_k^*T_{k}V_k$. We see that, although
$\hat{Q}_k$ is still column orthonormal, $Q^+_k=\hat{Q}_k\hat{V}_k$
is not as $\hat V_k$ is not orthogonal any longer when some rows are
deleted from the orthogonal matrix $V_k$. As a result,
$Q_m^+=Q_kV_k(:,1:m)$ is not column orthonormal, and
(\ref{resrt-3}) is not an $m$-step GSOAR decomposition any longer.
This means that implicit restarting fails to work whenever deflation
occurs.

In what follows we present
an effective approach to cure deflation so as to recover a standard GSOAR
decomposition, making implicit restarting always applicable unconditionally.

Note that $\hat{V}_k$ is a $(k-j)\times k$ of rank $k-j$.
Without loss of generality, we assume that
the first $k-j$ columns of $\hat{V}_k$ are linearly independent,
i.e., the matrix $\hat{V}_{k1}$ consisting the first $k-j$ columns of $\hat{V}_k$
is nonsingular. Write $\hat{V}_k=[\hat{V}_{k1},\hat{V}_{k2}]$. We
compute the QR decomposition of $\hat{V}_{k1}$ using the Matlab built-in function
${\sf qr.m}$ and obtain the decomposition of form
\begin{equation}\label{gm}
\hat{V}_k=U_kR_k=[U_{k-j},0]\left[\begin{array}{cc}
R_{k-j}&R_{12} \\
0& I
\end{array}
\right],
\end{equation}
where $\hat{V}_{k1}=U_{k-j}R_{k-j}$ is the QR decomposition of $\hat{V}_{k1}$ and
$R_{12}=U_{k-j}^*\hat{V}_{k2}$, and $I$ is the identity matrix of order $j$,
so that $R_k$ is nonsingular and upper triangular.

Noting that $U_k=\hat{V}_kR_k^{-1}$ and
right multiplying (\ref{defl-1}) by $R_k^{-1}$, we get
\begin{equation}
\left[ \begin{array}{cc}A&B \\I & 0\end{array}\right] \left[
\begin{array}{c}\hat{Q}_kU_k \\P_{k}V_kR_k^{-1}
\end{array}\right] =
\left[ \begin{array}{c}\hat{Q}_kU_k
\\P_{k}V_kR_k^{-1}\end{array}\right]
R_kT_{k}^+R_k^{-1} + t_{k+1k}\left [
\begin{array}{c} q_{k+1}\\p_{k+1}\end{array} \right]
e^*_{k}V_kR_k^{-1}. \label{defl-3}
\end{equation}
Since $R_k^{-1}$ is upper triangular, $R_kT_{k}^+R_k^{-1}$ is
Hessenberg. Note that $V_k$ has only $p=k-m$ nonzero subdiagonals.
Then the first possible nonzero entry $\tilde{\beta}$ of $e_k^*V_k$
is in position $m$ and
$$
t_{k+1k}e^*_{k}V_kR_k^{-1}=(0,\ldots,0,\tilde{\beta}, b^T)
$$
with $\tilde{\beta} = t_{k+1k}V_k(k,m)/e_m^*R_ke_m$.
Equating the first $m$ columns on two sides of (\ref{defl-3}),
we obtain
\begin{equation} \left[
\begin{array}{cc}A&B \\I & 0\end{array}\right] \left[
\begin{array}{c}\tilde{Q}_m^{+} \\\tilde{{P}}_m^{+}\end{array}\right]=
\left[ \begin{array}{c}\tilde{Q}_m^{+}
\\\tilde{P}_m^{+}\end{array}\right]\tilde{T}_m^{+}+\beta_m^+
\left [ \begin{array}{c} q_{m+1}^{+}\\p_{m+1}^{+}\end{array} \right]
e^*_m, \label{defl-4}
\end{equation}
where $\tilde{Q}_m^{+} = \hat{Q}_kU_k(:,1:m)$,
$\tilde{P}_m^+=P_kV_k(:,1:m)R_m^{-1}$ with $R_m$ the $m\times m$
leading principal matrix of $R_k$, $\tilde{T}_m^{+}$ the $m\times m$
leading principal matrix of $R_kT_k^+R_k^{-1}$, and
\begin{eqnarray}
\left [
\begin{array}{c} q_{m+1}^{+}\\p_{m+1}^{+}\end{array} \right]& =&
\frac{1}{\beta_m^+}f_m^+= \tilde{t}_{m+1m}^{+}\left[
\begin{array}{c}\hat{Q}_kU_k \\P_{k}V_kR_k^{-1}
\end{array}\right]e_{m+1}+ \tilde{\beta}
\left [ \begin{array}{c}
q_{k+1}\\p_{k+1}\end{array} \right],\\
\beta_m^+&=&\|\tilde{t}_{m+1m}^{+}\hat{Q}_kU_ke_{m+1}+ \tilde{\beta}q_{k+1} \|.
\end{eqnarray}

(\ref{gm}) indicates that the column orthonormality of $U_k(:,1:m)$ is guaranteed
whenever $m\leq k-j$, i.e., $j\leq k-m$. This means that
$\tilde{Q}_m^+=\hat{Q}_kU_k(:,1:m)$ is column orthonormal,
provided that the number $j$ of deflations during the last cycle of
GSOAR procedure does not exceed $k-m$. If $m>k-j$, the first $k-j$ columns of
$\tilde Q^+_m$ are orthonormal and the last $m-(k-j)$ columns of
$U_k$ are zero, so that the last $m-(k-j)$ columns of $\tilde Q^+_m$
are zero. As a result, there are $m-(k-j)$ deflations in (\ref{defl-4}).
For either $m\leq k-j$ or $m>k-j$, it is trivial to justify that
$(\tilde{Q}_m^+)^*q_{m+1}^+=0$. Therefore, by curing deflations in
implicit restarts, we have obtained a truly $m$-step GSOAR
decomposition (\ref{defl-4}).

\section{Numerical experiments}\label{s5}

In this section we report numerical examples to illustrate the
practicability of IGSOAR and IRGSOAR and the superiority of IRGSOAR to
IGSOAR. Meanwhile, we also compare them with the corresponding counterparts
IRSGA and IRRSGA proposed in \cite{wqhuang} for some test problems.
In addition, we compare IGSOAR and IRGSOAR with the Matlab function {\sf eigs},
the implicitly restarted Arnoldi method with exact shifts used, which is
directly applied to the linearization problem (\ref{5biao}).
All the experiments were run on Intel(R)Core(TM)i5-3470s CPU 2.9GHz,
RAM 4G using Matlab R2012b with $\epsilon_{\rm
mach}=2.22\times 10^{-16}$ under the Windows 7 system.

We list CPU timings (in second) of the three main parts
abbreviated as `SOAR', `SMALL' and `IMRE', where `SOAR' denotes
the CPU time of the first cycle of GSOAR procedure plus standard
extensions of the GSOAR decomposition from step $m+1$ to step $k$ for all
the other cycles, `SMALL' is the CPU time of forming the projected QEP,
solving them and computing residuals of approximate eigenpairs, and 'IMRE'
is the CPU time of performing all implicit QR iterations and generating
the $m$-step GSOAR decompositions for all cycles. In addition, we
use `restarts' and `CPU time' to denote the number of restarts and
the total CPU time of IGSOAR, IRGSOAR and {\sf eigs}, respectively.

For each example, we used the same starting vector generated randomly
in a uniform distribution for IGSOAR and IRGSOAR. We transformed the
projected QEP (\ref{05-pro}) to the generalized eigenvalue problem~(\ref{projgep})
and solved it by the QZ algorithm, i.e., the Matlab built-in function
{\sf eig.m}. We recovered an eigenvector $g$
of QEP (\ref{05-pro}) from either the first $k$ components or the
last $k$ components of $[\theta g^T,g^T]^T$. From the backward error
analysis \cite{higham}, it is preferable to take the first $k$ ones
if $|\theta|\geq 1$ and the last $k$ ones if $|\theta|<1$. We
adopted this choice.

For {\sf eigs}, we used the same $k$ as that in IGSOAR and IRGSOAR
to compute the same $m$ eigenpairs for each example. The CPU time of
{\sf eigs} did not include the time of computing
the LU decomposition of $M$, which is used when acting a matrix-vector
product in {\sf eigs} at each step.
The starting vector of {\sf eigs} was obtained by normalizing
$$
\left [ \begin{array}{c} q_1 \\ p_1 \end{array} \right ],
$$
where $q_1$ and $p_1$ were the vectors in Algorithm \ref{GSOAR}.
The number of shifts was the default value,
i.e., $p=k-(m+3)$. We also used $tol$ to denote
the stopping criterion used in {\sf eigs} for (\ref{5biao}).
Let $(\theta, y)$ be a converged eigenpair computed by {\sf eigs}, we set $y_1$ to be
the vector consisting of the first $n$ components of  $y$, and $y_2$ the vector
consisting of the last $n$ components of $y$. We then
computed the relative residual norms (\ref{stop}) of $(\theta,y_1)$ and $(\theta,y_2)$
and took the smaller one as the residual norm of {\sf eigs} for QEP (\ref{05-qep}).
'$Res_{\min}$ ' and  '$Res_{\max}$ ' recorded the minimum and maximum relative
residual norms (\ref{stop}) obtained in this way for all the converged eigenpairs
for (\ref{5biao}). The maximum number of restarts is limited to 50.

{\bf Example 1}. We consider the damped vibration mode of an acoustic
fluid confined in a cavity with absorbing walls capable of dissipating
acoustic energy \cite{huangjialin}. We take the same geometrical data as
in \cite{huangjialin}. The QEP is
$$
\lambda^2M_uu +(\alpha + \lambda \beta)A_u + K_uu=0,
$$
where $\alpha =5\times 10^4N/m^3$, $\beta = 200 Ns/m^3$, and the
order $n = 46548$.

By taking $tol=10^{-14}$ and two sets of parameters $k=30, p=7$
and $k=30,p=5$, we used IRGSOAR and IGSOAR to
compute the twenty eigenvalues nearest to the complex target $\sigma=25 + 18\pi i$
and the corresponding eigenvectors of the above QEP.
Table~\ref{5-T1} reports the results obtained, and Figure~\ref{figure5-1}
describes the convergence processes of two algorithms, depicting the maxima
of relative residual norms of $m$ approximate eigenpairs versus restarts.

We see from Table~\ref{5-T1} and Figure~\ref{figure5-1} that
two algorithms were efficient. However, as far as both restarts and CPU timings
are concerned, IRGSOAR was twice as fast as IGSOAR for $k=30$ and $p=5$, and the
former was also considerably faster than the latter for $k=30$ and $p=7$.
Furthermore,
we observe from the figure that the residual norm of IRGSOAR was smaller than
that of IGSOAR substantially at each cycle, indicating
that the refined Ritz vectors can be considerably more accurate than
the Ritz vectors. We find that for the same $k$, the value of $p$ has an effect
on the overall performance of IGSOAR and IRGOAR. For this example, we took two $p$
smaller than $k-m=10$. It is seen that the effect is
marginal for IRGSOAR, while it is relatively essential for IGSOAR.
In addition, we remark that the most consuming cost was paid to the SOAR procedure,
but the explicit computation and solutions of
all small QEP also occupied quite portion of the CPU time cost. The CPU time
'IMRE' of implicit restarting consumed least but could not be negligible.

\begin{table}[!hbp]
\begin{center}
\caption{Example 1, $tol$=$10^{-14}$ }
\label{5-T1}
\begin{tabular}{|c|c|c|c|c|c|c|c|}\hline
{Algorithm }&  $k$ & $p$& restarts& {CPU time} &  {SOAR} & {SMALL} & IMRE
 \\\hline
IRGSOAR & 30 &7 & 3& 12.55 & 8.56   & 2.58 & 1.34 \\\hline
IGSOAR & 30  &7  & 5& 16.86 & 10.99 & 3.12 & 2.68 \\\hline
IRGSOAR & 30 &5 & 3& 11.96  & 7.93  & 2.58 & 1.31 \\\hline
IGSOAR & 30  &5  & 7& 21.01 & 12.06 & 4.57 & 4.27 \\\hline
\end{tabular}
\centerline{The results obtained by {\sf eigs}}
\begin{tabular}{|c|c|c|c|c|c|c|c|c|}\hline
{$tol$ }&  $k$ & CPU & restarts & $Res_{\min}$ & $Res_{\max}$ \\ \hline
$10^{-6}$ & 30 &27.89 &15  & $4.15\times 10^{-16}$ &$9.25\times 10^{-15}$  \\ \hline
$10^{-8}$ & 30 &--     &50  & -- &--  \\ \hline
\end{tabular}
\end{center}
\end{table}

\begin{figure}[!hbp]
\begin{center}
\includegraphics[height=5cm,width=6cm]{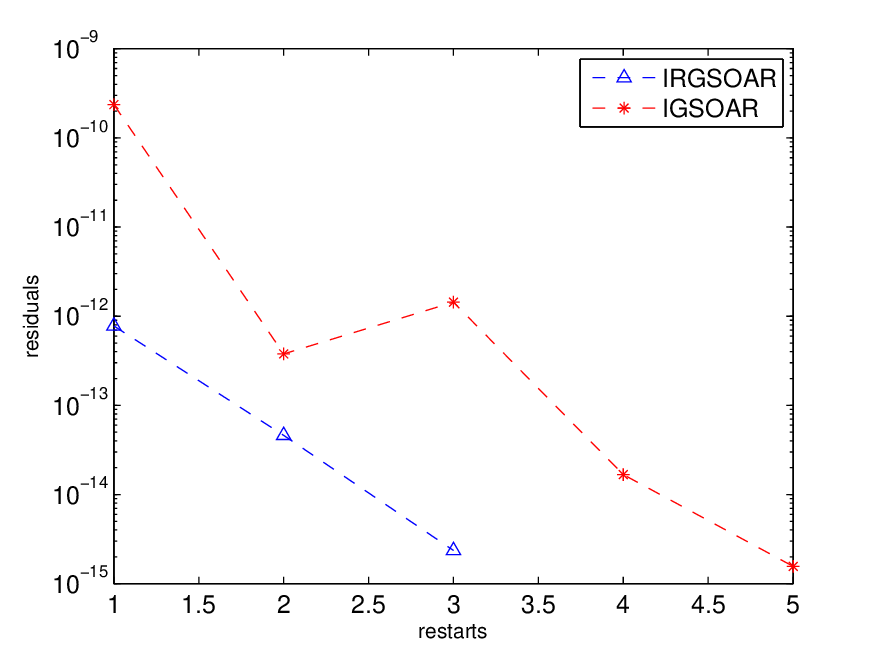}
\includegraphics[height=5cm,width=6cm]{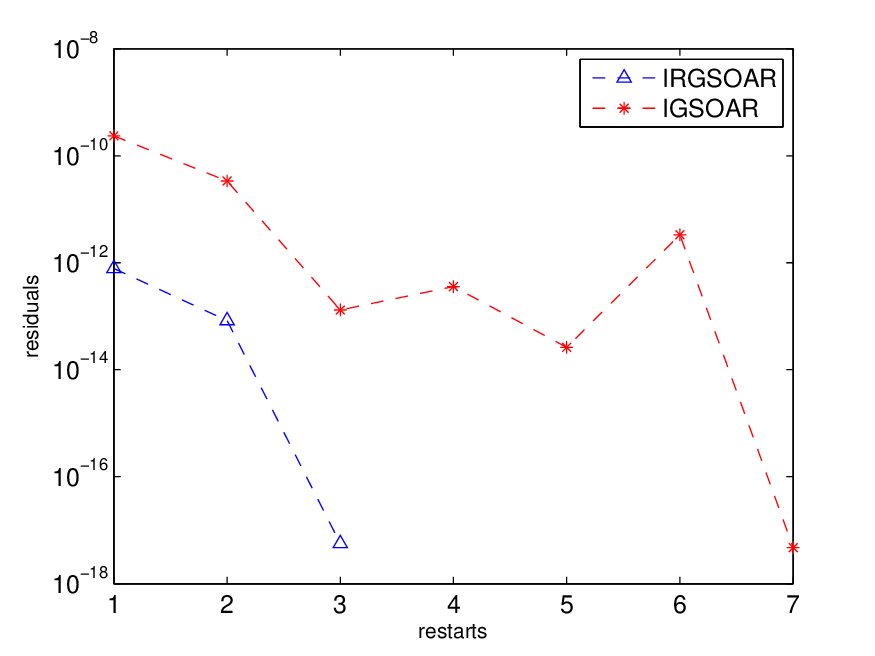}
\caption{Example 1. Residuals versus restarts.
Left: $k = 30,\ p= 7$; right: $ k = 30,\ p = 5$.}\label{figure5-1}
\end{center}
\end{figure}

For this example, by taking $tol=10^{-6}$, we found that
{\sf eigs} was much more costly than IGSOAR and IRGSOAR
to converge and the approximate eigenpairs were as accurate as
those obtained by the latter two algorithms,
while, for $tol=10^{-8}$, it failed to converge after 50 restarts.
We should point out that our codes are programmed
in the Matlab language and may not be optimized while {\sf eigs}
is programmed in C language and optimized. This means that for the same
$k$ each restart of {\sf eigs} should be more time consuming than
that of IGSOAR and IRGSOAR since {\sf eigs} is much more
expensive than IGSOAR and IRSOAR in the orthogonalization
of Arnoldi vectors. As a result, in all the experiments
the number of restarts is
more reasonable to compare the computational efficiency of
these three algorithms. It is worthwhile to mention that for this example
a relatively big $tol=10^{-6}$ for {\sf eigs} delivered very accurate
eigenpairs of QEP (\ref{05-qep}) and a smaller $tol$ is unnecessary.

{\bf Example 2}. This problem and arises in a model of the
concrete structure supporting a machine assembly \cite{Betcke,
Feriani} and has the form $Q(\lambda)x =(\lambda^2M + \lambda C + (1 + i\mu)K)x=0$.
The matrices are of order 2472, where $M$ is real diagonal,
$C$, the viscous damping matrix, is pure imaginary and
diagonal, $K$ is complex symmetric, and the factor $1 + i\mu$ adds
uniform hysteretic damping. We use the command {\sf nlevp(`concrete', 0.04)}
in \cite{Betcke} to generate the complex symmetric coefficient matrices.
Thus problem was tested in \cite{wqhuang}.

We ran IRGSOAR and IGSOAR to compute the ten eigenvalues nearest to the origin
by taking $tol= 10^{-14}$ and the same $k = 20$, two $p=7$ and $5$. Table~\ref{4-T1}
and Figure~\ref{figure4-1} reported the results, from which it can be seen that two
algorithms worked very well and IRGSOAR was a little more efficient than IGSOAR
in terms of both restarts and CPU timings. We remark that, for this problem,
the corresponding algorithms IRSGA and IRRSGA in \cite{wqhuang} both used four restarts
to achieve the convergence for the same $k=20$ and $tol$.
Note that they use the F-norm in the denominator of (\ref{stop}), which means that
for the same $tol$ our convergence tolerance is smaller. Therefore, for $p=7$,
IGSOAR was (at least) as efficient as IRSGA, and IRGSOAR
was faster than IRRSGA. For $p=5$, IGSOAR used five restarts for a smaller
stopping tolerance than that used by IRSGA, and IRGSOAR used four restarts. This
demonstrates that, for this problem, IGSOAR and IRGSOAR were as efficient as IRSGA and
IRRSGA, respectively.  It is clear that two different $p$ affected the overall
efficiency of each algorithm only marginally. Finally, we observe that,
unlike Example 1, the main cost of each algorithm
was paid to the GSOAR procedure and overwhelmed "SMALL" and "IMRE".

\begin{table}[h]
\begin{center}
\caption{Example 2, $tol$=$10^{-14}$ }
\label{4-T1}
\begin{tabular}{|c|c|c|c|c|c|c|c|c|c|c|c|}\hline
{Algorithm }&  $k$ & $p$& restarts& {CPU time} &  {SOAR} & SMALL& IMRE \\\hline
  IRGSOAR & 20 & 7  & 3& 0.78 & 0.62 & 0.10 & 0.03  \\\hline
  IGSOAR  & 20 & 7  & 4 & 1.03 & 0.91 & 0.07 & 0.04  \\\hline
  IRGSOAR & 20 & 5  & 4 &0.90 & 0.72 & 0.12 & 0.05  \\\hline
  IGSOAR  & 20 & 5  & 5& 0.93 & 0.77 & 0.09 & 0.07  \\\hline
\end{tabular}\\
\centerline{The results obtained by {\sf eigs}}
\begin{tabular}{|c|c|c|c|c|c|c|c|c|}\hline
{$tol$ }&  $k$ & CPU time & restarts &  $Res_{\min}$  & $Res_{\max}$  \\ \hline
$10^{-8}$  & 20 &1.03 &7  & $1.95\times10^{-18}$ &  $1.26\times10^{-14}$   \\ \hline
$10^{-10}$ & 20 &1.25 &9  & $1.88\times10^{-18}$ &  $1.18\times10^{-16}$ \\ \hline
$10^{-14}$ & 20 &1.70 &12 & $1.35\times10^{-18}$ &  $5.15\times10^{-18}$  \\ \hline
\end{tabular}
\end{center}
\end{table}
\begin{figure}[!h]
\begin{center}
\includegraphics[height=5cm,width=6cm]{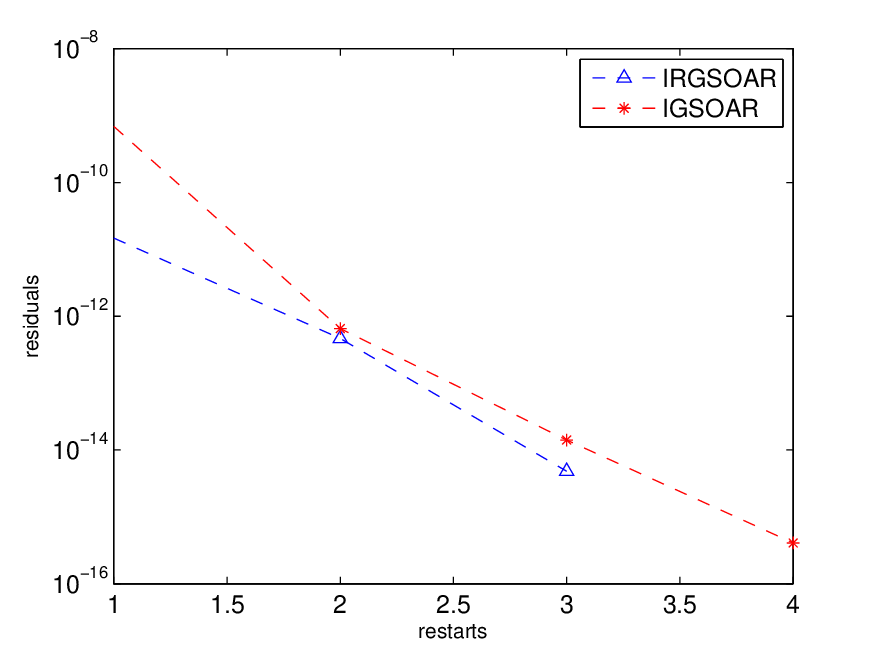}
\includegraphics[height=5cm,width=6cm]{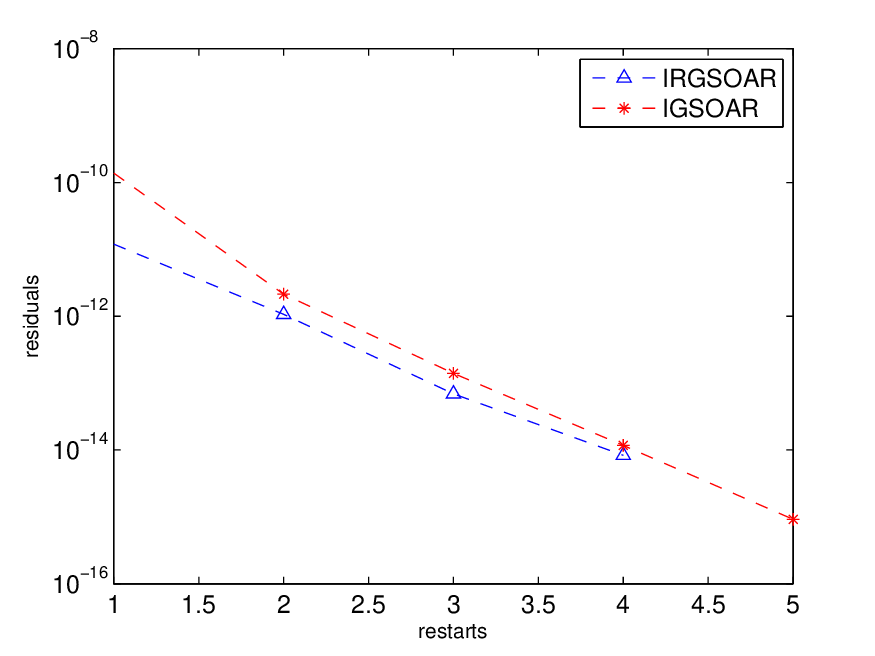}
\caption{Example 2. Residuals versus restarts;
 Left: $k = 20,\ p= 7$; right: $ k = 20,\ p = 5$.}\label{figure4-1}
\end{center}
\end{figure}

We also report the results obtained by {\sf eigs} for three $tol$  and list them in
Table~\ref{4-T1}. It is seen from Table~\ref{4-T1} that both IGSOAR and IRGSOAR
performed much better than {\sf eigs},
and they used much less CPU time and fewer restarts to compute the desired
eigenpairs with much higher accuracy. As $Res_{\min}$ and $Res_{\max}$ indicated, 
the accuracy of the converged eigenpairs obtained by {\sf eigs} with three greatly 
varying $tol$ essentially had no difference as the approximate eigenpairs of 
QEP (\ref{05-qep}), and their relative residual norms were already at the level of 
$\epsilon_{\rm mach}$ for $tol=10^{-8}$.

{\bf Example 3}. This example is from \cite{Betcke} and tested in \cite{wqhuang}
(cf. Example 6.3). We tested IRGSOAR and IGSOAR for the following
cases (a) and (b) by taking $tol=10^{-14}$.

Case (a): Acoustic 1D.
This example arises from the finite element discretization of the time harmonic
wave equation $-\triangle p-(2\pi f/c)^2 p=0$. Here, $p$ denotes the pressure, $f$
is the frequency, $c$ is the speed of sound in the medium, and $\xi$
is the (possibly complex) impedance. On the
domain $[0, 1]$ with $c=1$, the $n\times n$ matrices $M$, $D$, and $K$ are defined by
$$
M=-4\pi^2\frac{1}{n}\left(I-e_ne_n^T\right),\
D=2\pi i\frac{1}{\xi}e_ne_n^T,\ K=n\left({\rm tridiag}(-1,2,-1)-e_ne_n^T\right).
$$
We use {\sf nlevp(¡®acoustic\_wave\_1d¡¯,5000,1)} to generate matrices $M,\ D$
and $K$ with size $n=5000$.

Just as in \cite{wqhuang},
we computed the six eigenvalues nearest to the origin with $k=12, p=5$ and 3.
Table~\ref{5-a} reports the results, and Figure \ref{figure7-1} depicts
the convergence processes of two algorithms. From the figure we see that,
for the same $k$ and two $p$, IRGSOAR
and IGSOAR used two and three cycles, respectively. As indicated in
\cite{wqhuang}, for the same $k=12$ and a little larger convergence tolerance,
the corresponding implicitly restarted algorithms IRSGA and IRRSGA
both used three cycles. So IRGSOAR was a little better than IRRSGA.
Regarding CPU time, since $M, D$ and $K$ are very sparse, the
CPU timings of the GSOAR procedure and implicit restarting are comparable,
and less than `SMALL'.

\begin{table}[h]
\begin{center}
\caption{Example 3(a), $tol$=$10^{-14}$ }
\label{5-a}
\begin{tabular}{|c|c|c|c|c|c|c|c|}\hline
{Algorithm }&  $k$ & $p$& restarts & {CPU time} &  {SOAR} & {SMALL} & IMRE
 \\\hline
IRGSOAR & 12 &5 & 2 & 0.12 & 0.02 & 0.06 & 0.02 \\\hline
IGSOAR  & 12 &5 & 3 & 0.12 & 0.02 & 0.05 &  0.03 \\\hline
IRGSOAR & 12 &3 & 2 & 0.10 & 0.02 & 0.06 & 0.01 \\\hline
IGSOAR  & 12 &3 & 3 & 0.10 & 0.02 & 0.05 &  0.02 \\\hline
\end{tabular}
\centerline{The results obtained by {\sf eigs}}
\begin{tabular}{|c|c|c|c|c|c|c|c|c|}\hline
{$tol$ }&  $k$ & CPU time & restarts &  $Res_{\min}$  & $Res_{\max}$  \\ \hline
$10^{-8}$  & 12 &0.27 &10 & $0.32\times10^{-18}$ &  $0.97\times10^{-13}$   \\ \hline
$10^{-10}$ & 12 &0.72 &22 & $0.40\times10^{-18}$ &  $0.20\times10^{-15}$ \\ \hline
$10^{-14}$ & 12 &0.75 &23 & $0.81\times10^{-18}$ &  $0.83\times10^{-18}$ \\ \hline
\end{tabular}
\end{center}
\end{table}

\begin{figure}[!hbp]
\begin{center}
\includegraphics[height=5cm,width=6cm]{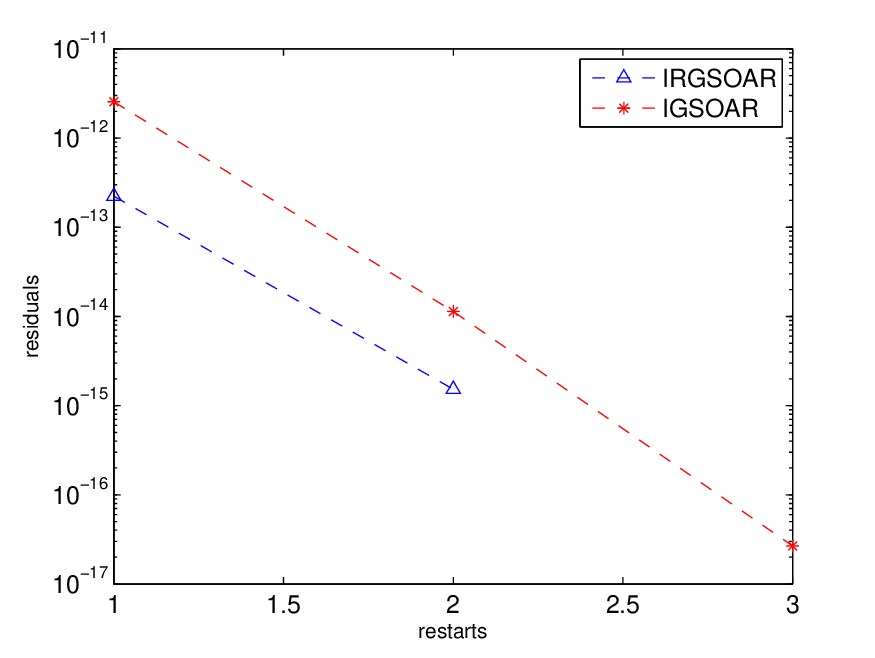}
\includegraphics[height=5cm,width=6cm]{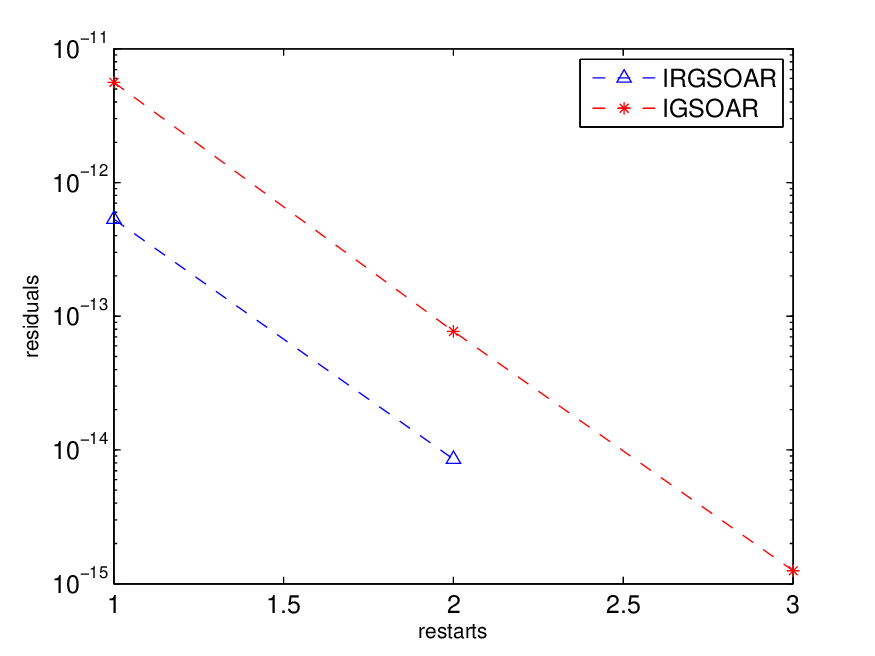}
\caption{ Example 3(a).
 Left: $k = 12,\ p= 5$; right: $ k = 12,\ p = 3$.}  \label{figure7-1}
\end{center}
\end{figure}

As we have seen, the eigenpairs obtained by {\sf eigs} had similar accuracy to
those obtained by IGSOAR and IRGSOAR with three greatly varying $tol$,
and all of them converged to the level of machine precision.  However,
Table~\ref{5-a} clearly shows that IGSOAR and IRGSOAR were much more efficient than
{\sf eigs}.

Case (b): Acoustic 2D. This example is a two-dimensional acoustic wave equation
on $[0,1]\times [0,1]$. The coefficient matrices $M,\ D$ and $K$ are given by
\begin{eqnarray*}
M&=&-4\pi^2 h^2 I_{q-1}\otimes \left(I_q-\frac{1}{2}e_qe_q^T\right),\
D=2\pi i\frac{h}{\xi} I_{q-1}\otimes (e_qe_q^T),\\
K&=&I_{q-1}\otimes D_q+T_{q-1}\otimes \left(-I_q+\frac{1}{2}e_qe_q^T
\right).
\end{eqnarray*}
where $h$ denotes the mesh size, $q=1/h$, $\otimes$ denotes the Kronecker product,
$\xi$ is the (possibly complex) impedance, $D_q={\rm tridiag}(-1, 4,-1)- 2e_qe_q^T$,
and $T_{q-1}={\rm tridiag}(1, 0, 1)$. We use
{\sf nlevp(`acoustic\_wave\_2d',90,0.1$\times$
1i)} to get the real symmetric matrices $(M,D,K)$. The matrix size is given
by $n=8010$.

As in \cite{wqhuang}, we computed  the six eigenvalues nearest to the origin
with $k=12, p=5$ and 3. Table~\ref{5-b} and Figure~\ref{figure7-2} give
the results and convergence processes of two algorithms, respectively. It
is seen that IRGSOAR and IGSOAR used seven and eleven restarts for $p=5$,
respectively, and both of them used ten cycles for $p=3$.
Therefore, two algorithms were efficient,
and IRGSOAR could be more efficient than IGSOAR. We find that
both IGSOAR and IRGSOAR were more efficient than IRSGA
and IRRSGA \cite{wqhuang}, where the latter ones used eleven and
twelve cycles, respectively.

\begin{table}[h]
\begin{center}
\caption{Example 3(b), $tol$=$10^{-14}$ }
\label{5-b}
\begin{tabular}{|c|c|c|c|c|c|c|c|}\hline
{Algorithm }&  $k$ & $p$& restarts &{CPU time} &  {SOAR} & {SMALL} & IMRE
 \\\hline
IRGSOAR & 12 &5 & 7 & 1.21 & 0.59  & 0.39 & 0.19 \\\hline
IGSOAR  & 12 &5 & 11 & 1.57 & 0.85  & 0.36 & 0.33 \\\hline
IRGSOAR & 12 &3 & 10 & 1.18 & 0.45  & 0.45 & 0.24 \\\hline
IGSOAR  & 12 &3 & 10 & 0.94 & 0.46  & 0.25 & 0.22 \\\hline
\end{tabular}

\centerline{The results obtained by {\sf eigs}}
\begin{tabular}{|c|c|c|c|c|c|c|c|c|}\hline
{$tol$ }&  $k$ & CPU time & restarts &  $Res_{\min}$  & $Res_{\max}$  \\ \hline
$10^{-8}$  & 12 &0.84 &8 & $0.38\times10^{-17}$ &  $0.48\times10^{-13}$   \\ \hline
$10^{-10}$ & 12 &0.83 &8 & $0.40\times10^{-17}$ &  $0.14\times10^{-13}$ \\ \hline
$10^{-14}$ & 12 &0.99 &11 & $0.32\times10^{-17}$ &  $0.53\times10^{-16}$ \\ \hline
\end{tabular}
\end{center}
\end{table}

\begin{figure}[!hbp]
\begin{center}
\includegraphics[height=5cm,width=6cm]{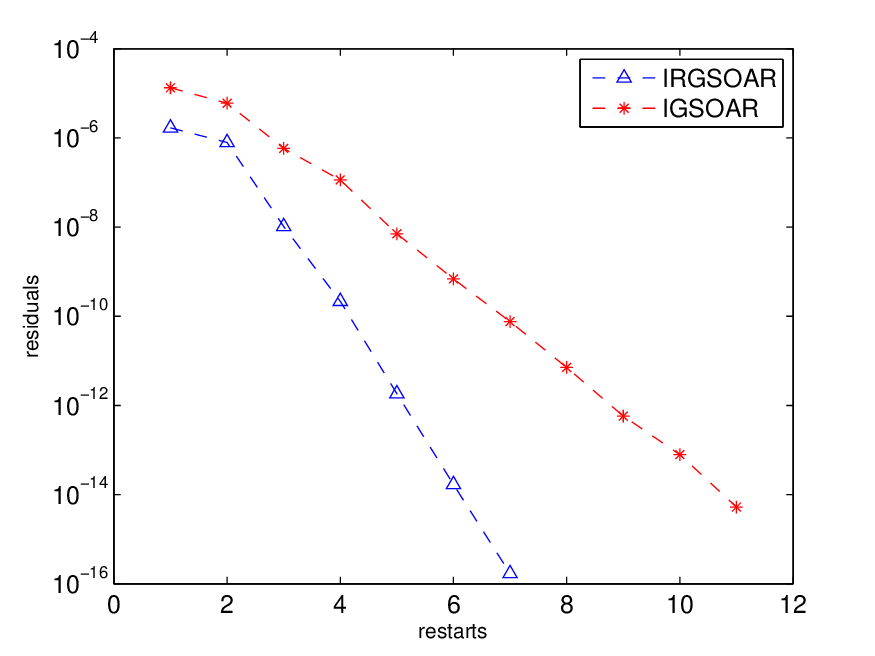}
\includegraphics[height=5cm,width=6cm]{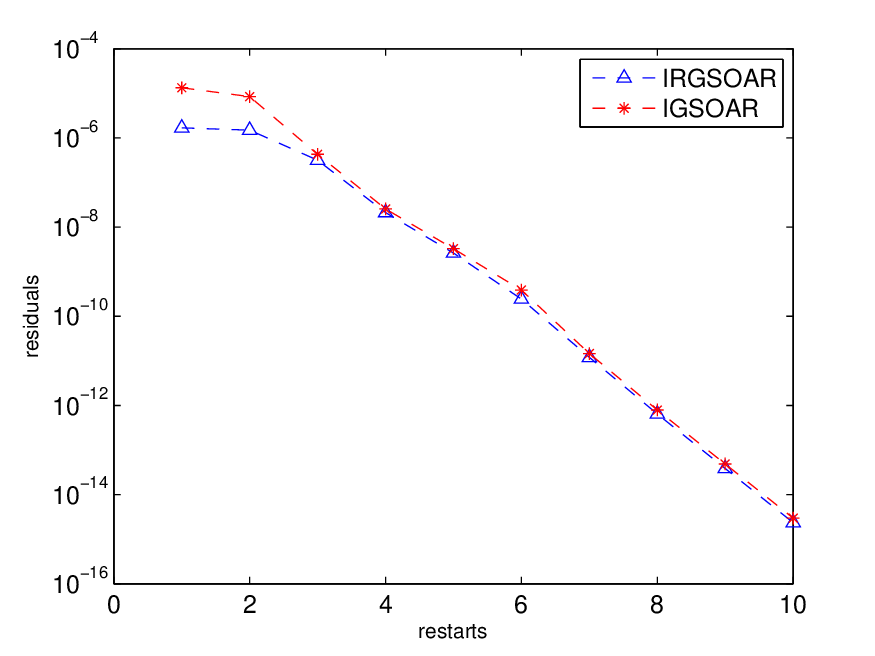}
\caption{Example 3(b).
 Left: $k = 12,\ p= 5$; right: $ k = 12,\ p = 3$.}  \label{figure7-2}
\end{center}
\end{figure}

For this problem, unlike the previous examples, {\sf eigs}
performed very well and was comparable to IGSOAR and IRGSOAR
in terms of the accuracy and the computational efficiency.

{\bf Example 4}. This QEP arises in an
$n$-degree-of-freedom damped mass-spring system \cite{tm}. By taking
$m_i=1$ and letting all the springs (respectively, dampers) have the
same constant $\kappa$ (respectively, $\tau$) except
$\kappa_1=\kappa_n=2\kappa$ and $\tau_1=\tau_n=2\tau$, the resulting
matrices are
$$
M=I,~~C=\tau\cdot {\rm tridiag}(-1,3,-1),~~
K=\kappa\cdot {\rm tridiag}(-1,3,-1),
$$
which are very sparse. We took $n=5000,\kappa=5$ and $\tau=10$ and
were interested in the six eigenvalues nearest to the complex target
$\sigma= -13 + 0.4i$ and the corresponding eigenvectors.

For $tol=10^{-10}$, we tested IRGSOAR and IGSOAR for $k=40, p=23$ and 28.
Table~\ref{2-T1} lists the results, and
Figure~\ref{figure2-1} depicts the convergence processes
for two sets of parameters $k$ and $p$.

\begin{table}[h]
\begin{center}
\caption{Example 4, $tol$=$10^{-10}$}
\label{2-T1}
\begin{tabular}{|c|c|c|c|c|c|c|c|c|c|c|c|}\hline
{Algorithm }&  $k$ & $p$& restarts &{CPU time} &  {SOAR} & SMALL& IMRE
\\\hline
IRGSOAR & 40 &23 &41 & 11.01 & 1.78 & 4.74 & 4.39 \\\hline
IGSOAR  & 40 &23 &44 & 9.02  & 1.88 & 2.33 & 4.73  \\\hline
IRGSOAR & 40 &28 &39 & 9.96  & 1.79 & 4.07 & 4.01  \\\hline
IGSOAR  & 40 &28 &47 & 9.50  & 2.18 & 2.47 & 4.78 \\\hline
\end{tabular}
\centerline{The results obtained by {\sf eigs}}
\begin{tabular}{|c|c|c|c|c|c|c|c|c|}\hline
{$tol$ }&  $k$ & CPU time & restarts &  $Res_{\min}$  & $Res_{\max}$  \\ \hline
$10^{-8}$  & 40 &6.05 &31 & $0.37\times10^{-3}$ &  $0.37\times10^{-3}$ \\ \hline
$10^{-10}$ & 40 &9.36 &47 & $0.37\times10^{-3}$  &  $0.37\times10^{-3}$ \\ \hline
\end{tabular}

\end{center}
\end{table}
It can be found from Table~\ref{2-T1} and Figure~\ref{figure2-1} that two algorithms
worked quite well. Compared with Examples 1--3, much more restarts were needed now;
for the given $k$, two different $p$ did not make much difference on restarts and
CPU timings of two algorithms. Furthermore, IRGSOAR and IGSOAR are similarly
efficient, and the former used a little fewer restarts but more CPU time than IGSOAR.
Since the matrices in this
QEP are very sparse, it appears that performing the SOAR procedure in each algorithm
was not dominant, and instead it was considerably less costly than the explicit
computation and solutions of all small QEP
and implicit restarting, as indicated by Table~\ref{2-T1}.

\begin{figure}[]
\begin{centering}
\includegraphics[height=5cm,width=6cm]{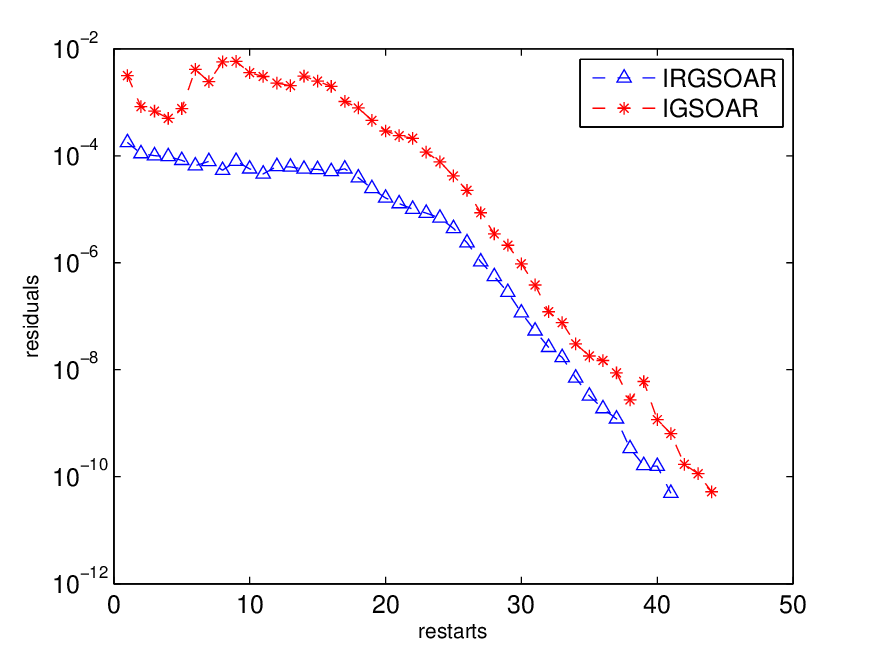}
\includegraphics[height=5cm,width=6cm]{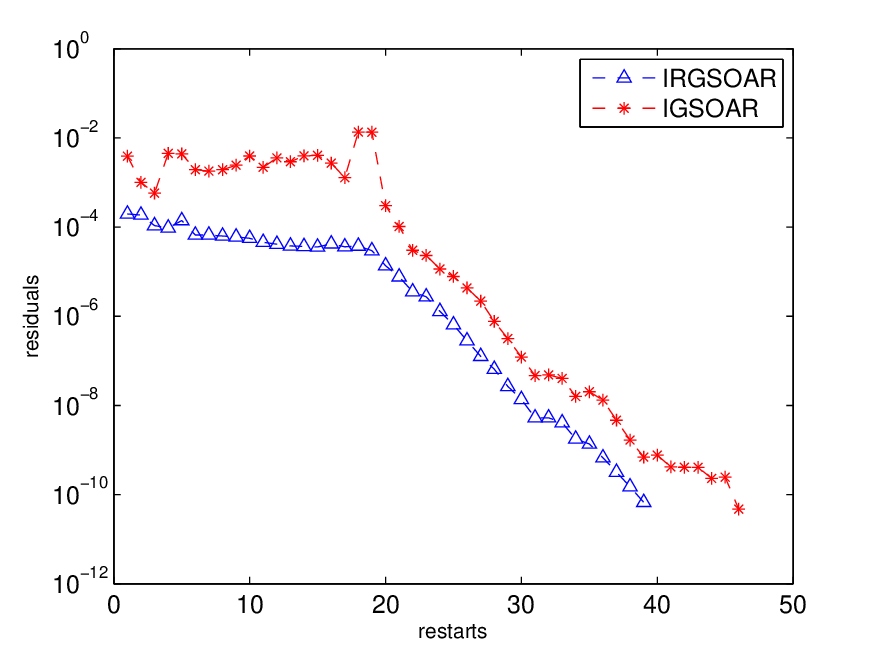}
\caption{Example 4. Residuals versus restarts. Left: $k = 40,\ p
= 23$; right: $ k = 40,\ p = 28$.}  \label{figure2-1}
\end{centering}
\end{figure}

For this example, unlike all the previous examples,
for given two $tol$ similar to that used by
IGSOAR and IRGSOAR, {\sf eigs} used comparable restarts and the CPU time
to IGSOAR and IRGSOAR, but it computed the desired eigenpairs with much poorer
accuracy, so, as a whole, it is considerably inferior to IGSOAR and IRGSOAR.
An important observation is that improving the accuracy
of approximate eigenpairs of (\ref{5biao}) may be helpless to improve
their accuracy as the approximate eigenpairs of (\ref{05-qep}).
A comparison of this example and Example 1 reveals a remarkable
difference: {\sf eigs} with big $tol$ computed the desired eigenpairs
with the accuracy at the level of $\epsilon_{\rm mach}$ for Example 1,
while it with smaller $tol$ got the desired eigenpairs with much poorer
accuracy. So it is uncertain for us to choose a suitable $tol$ for
{\sf eigs} to compute the desired eigenpairs with a prescribed accuracy
in the sense of the stopping criterion (\ref{stop}) for QEP (\ref{05-qep}).

{\bf Example 5}. This problem comes from \cite{Betcke}.
It is a nonlinear eigenvalue problem modeling a
radio-frequency gun cavity that is of the form
$$T(\lambda)x = [K - \lambda M + i(\lambda- \sigma_1^2)^{1/2}W_1 +
i(\lambda- \sigma_2^2)^{1/2}W_2]x = 0,$$ where $M,\ K,\ W_1,\ W_2$ are
real symmetric matrices of size $9956 \times 9956$. From these
matrices, we constructed a QEP of the form
$$
(\lambda^2W_2 + \lambda M +K)x=0,
$$
which is purely for our test purpose. We used IRGSOAR and IGSOAR to compute
the six eigenvalues nearest to $\sigma=0.5+0.5i$ and the associated eigenvectors.
Table~\ref{3-T1} and Figure~\ref{figure3-1} reported the results.

\begin{table}[!h]
\begin{center}
\caption{Example 5, $tol$=$10^{-10}$ }
\label{3-T1}
\begin{tabular}{|c|c|c|c|c|c|c|c|c|c|c|c|c|}\hline
{Algorithm }&  $k$ & $p$& restarts & {CPU time} &  {SOAR} & SMALL& IMRE
\\\hline
IRGSOAR & 20 &5  &1 & 1.73  & 1.59 & 0.11 & 0.00  \\\hline
IGSOAR  & 20 &5  &9 & 6.10  & 4.93 & 0.58 & 0.54  \\\hline
IRGSOAR & 20 &11 &1 & 1.73  & 1.61 & 0.11 & 0.00 \\\hline
IGSOAR  & 20 &11 &3 & 3.77  & 3.42 & 0.19 & 0.13 \\\hline
\end{tabular}
\centerline{The results obtained by {\sf eigs}}
\begin{tabular}{|c|c|c|c|c|c|c|c|c|}\hline
{$tol$ }&  $k$ & CPU time & restarts &  $Res_{\min}$  & $Res_{\max}$  \\ \hline
$10^{-6}$ & 20 &2.73  &2 & $0.26\times10^{-10}$ &  $0.25\times10^{-7}$ \\ \hline
$10^{-8}$ & 20 &59.31 &50 & $0.25\times10^{-7}$  &  $0.21\times10^{-3}$ \\ \hline
\end{tabular}
\end{center}
\end{table}
For this example, two algorithms worked well. However, IRGSOAR exhibited
the very considerable superiority to
IGSOAR. We find the desired eigenpairs without restarting the algorithm
for given two sets of parameters $k$ and $p$ while IGSOAR
used nine and three cycles, respectively. In terms of CPU timings, IRGSOAR
was also a few times faster than IGSOAR. Furthermore, for this example, the CPU time
of the SOAR procedure dominated the CPU time cost of each algorithm.
On contrary to Example 4, for the given $k$, the smaller $p=5$ made IGSOAR use
considerably more restarts and CPU time, meaning that the choice of $p$ may
have considerable effects on
the overall performance of IGSOAR. However, this example and
Examples 1--4 illustrate that the effects of $p$ must be problem dependent,
and it is impossible to design a definite and general
effective way to select it.

\begin{figure}[h]
\begin{centering}
\includegraphics[height=5cm,width=6cm]{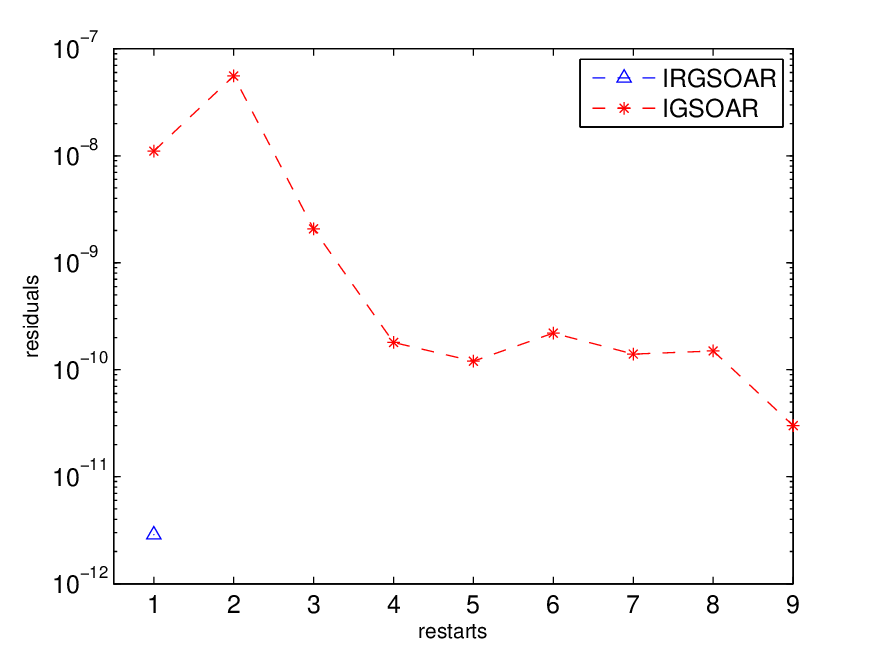}
\includegraphics[height=5cm,width=6cm]{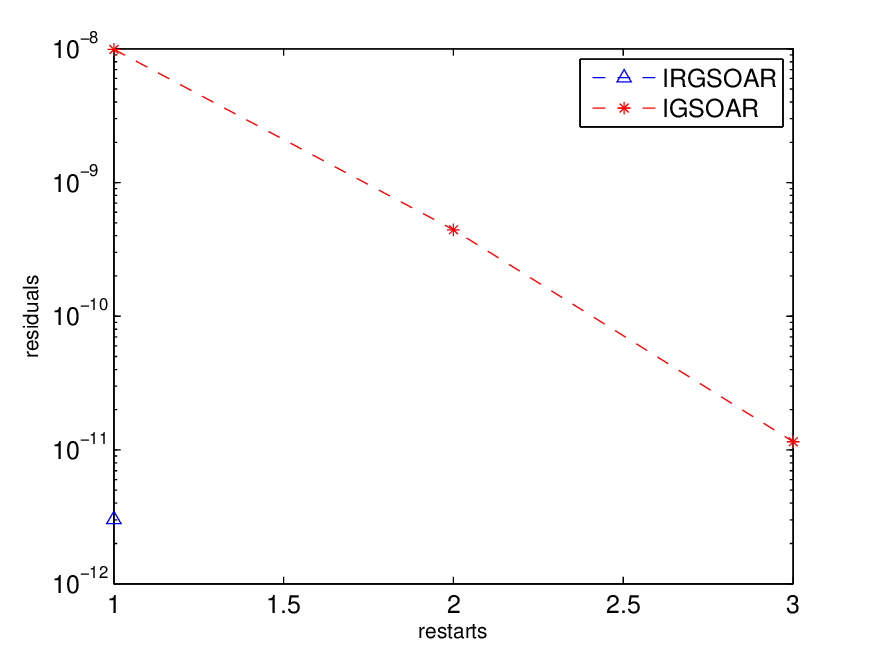}
\caption{Example 5. Residuals versus restarts. Left: $ k = 20,\ p = 5$;
right: $ k = 20,\ p = 11$.}\label{figure3-1}
\end{centering}
\end{figure}

In contrast, {\sf eigs} behaved not good for this example, and it used much more restarts
to achieve the convergence for the not much smaller $tol=10^{-8}$ than $10^{-6}$.
However, as approximate eigenpairs of QEP (\ref{05-qep}),
the converged eigenpairs with $tol=10^{-6}$ were substantially more accurate
than those with $tol=10^{-8}$. This is really bad because it shows that,
on the contrary to our common acceptance, that considerably more accurate
eigenpairs for the linearization
problem (\ref{5biao}) are not necessarily more accurate too for QEP (\ref{05-qep}).
This, together with Example 1 and Example 4, demonstrates that solving the
linearization problem (\ref{5biao}) directly
has serious uncertainty, as far as the accuracy is concerned.

\section{Conclusion}

We have considered generalized second-order Arnoldi method and its
refined version for solving the large QEP. The methods are structure-preserving
and applied to the QEP directly after
an orthonormal basis of the generalized second-order Krylov subspace
is generated by the GSOAR procedure.
To be practical, we have developed implicitly restarted algorithms
with certain exact and refined shifts proposed for two methods, respectively.
We have presented an efficient and reliable algorithm for computing the shift
candidates. Unlike Arnoldi type algorithms for the linear eigenvalue problem,
where the number of shift candidates are just that of shifts, for the QEP
the shift candidates are more than the shifts. we have discussed in detail how
to seek and determine reasonable shifts for each method. Also,
deflation may occur in the algorithms for the QEP, for which implicit
restarting is not applicable. To overcome this deficiency,
we have proposed an effective approach to cure deflation in implicit restarts,
so that implicit restarting can be used to the GSOAR procedure unconditionally.
We have tested our algorithms on a number of real-world problems. Numerical
experiments have demonstrated that two algorithms work well and the refined
algorithm can outperform the standard counterpart considerably. They also show
that our algorithms generally perform much better than {\sf eigs} in terms of the
accuracy or the computational efficiency.

\end{document}